\documentclass[12pt]{article}

\advance\voffset by -0.25cm

\usepackage{amsfonts,amsmath,amssymb,amsthm,bbm,eurosym,fancyhdr,graphicx,setspace,url}

\newtheorem{theorem}{Theorem}[section]

\newtheorem{definition}[theorem]{Definition}

\begin{document}

\newcommand{\real}{{\ensuremath{\mathbb{R}}}}
\newcommand{\dd}{{\ensuremath{\rm d}}}

\newcommand{\bS}{\mbox{\bf S}}

\newcommand{\cA}{{\ensuremath{\cal A}}}
\newcommand{\cD}{{\ensuremath{\cal D}}}
\newcommand{\cF}{{\ensuremath{\cal F}}}
\newcommand{\cH}{{\ensuremath{\cal H}}}
\newcommand{\cN}{{\ensuremath{\cal N}}}
\newcommand{\cP}{{\ensuremath{\cal P}}}

\newcommand{\hsp}{\hspace{0.2mm}}
\newcommand{\vsp}{\vspace{-2mm}}

\newcommand{\myA}{\mbox{\rm A}}
\newcommand{\myD}{\mbox{\rm D}}
\newcommand{\myI}{\mbox{\rm I}}
\newcommand{\myL}{\mbox{\rm L}}
\newcommand{\myO}{\mbox{\rm O}}
\newcommand{\myS}{\mbox{\rm S}}
\newcommand{\myT}{\mbox{\rm T}}
\newcommand{\myU}{\mbox{\rm U}}
\newcommand{\myV}{\mbox{\rm V}}

\newcommand{\myE}{{\ensuremath{\mathbb{E}}}}
\newcommand{\one}{\ensuremath{\mathbbm{1}}}

\newcommand{\Ers}{\frac{\myE_F \hsp [r(Y)]}{\myE_F \hsp [s(Y)]}}
\newcommand{\done}{\hfill $\Box$} 

\newcommand{\bysame}{\raisebox{1.5mm}{\underline{\hspace{2.75em}}}}
\newcommand{\marginal}[1]{\marginpar{\raggedright\scriptsize #1}}

\begin{center}
{

\LARGE \bf 

Making and Evaluating Point Forecasts

} 

\bigskip
{\bf \large Tilmann Gneiting}

\bigskip
{\bf Institut f\"ur Angewandte Mathematik \\ Universit\"at Heidelberg}

\bigskip
{\bf \today}
\end{center}

\vspace{2mm}
\begin{abstract}
Single-valued point forecasts continue to be issued and used in almost
all realms of science and society.  Typically, competing point
forecasters or forecasting procedures are compared and assessed by
means of an error measure or scoring function, such as the absolute
error or the squared error, that depends both on the point forecast
and the realizing observation.  The individual scores are then
averaged over forecast cases, to result in a summary measure of the
predictive performance, such as the mean absolute error or the (root)
mean squared error.  I demonstrate that this common practice can lead
to grossly misguided inferences, unless the scoring function and the
forecasting task are carefully matched.

Effective point forecasting requires that the scoring function be
specified a priori, or that the forecaster receives a directive in the
form of a statistical functional, such as the mean or a quantile of
the predictive distribution.  If the scoring function is specified a
priori, the forecaster can issue an optimal point forecast, namely,
the Bayes rule, which minimizes the expected loss under the
forecaster's predictive distribution.  If the forecaster receives a
directive in the form of a functional, it is critical that the scoring
function be consistent for it, in the sense that the expected score is
minimized when following the directive.  Any consistent scoring
function induces a proper scoring rule for probabilistic forecasts,
and a duality principle links Bayes rules and consistent scoring
functions.

A functional is elicitable if there exists a scoring function that is
strictly consistent for it.  Expectations, ratios of expectations and
quantiles are elicitable.  For example, a scoring function is
consistent for the mean functional if and only if it is a Bregman
function.  It is consistent for a quantile if and only if it is
generalized piecewise linear.  Similar characterizations apply to
ratios of expectations and to expectiles.  Weighted scoring functions
are consistent for functionals that adapt to the weighting in peculiar
ways.  Not all functionals are elicitable; for instance, conditional
value-at-risk is not, despite its popularity in quantitative finance.

\end{abstract}

\smallskip
\section{Introduction}  \label{sec:introduction}

In many aspects of human activity, a major desire is to make forecasts
for an uncertain future.  Consequently, forecasts ought to be
probabilistic in nature, taking the form of probability distributions
over future quantities or events (Dawid 1984; Gneiting 2008a).  Still,
many practical situations require single-valued point forecasts, for
reasons of decision making, market mechanisms, reporting requirements,
communications, or tradition, among others.

\subsection{Using scoring functions to evaluate point forecasts}  \label{sec:scoring}

In this type of situation, competing point forecasters or forecasting
procedures are compared and assessed by means of an error measure,
such as the absolute error or the squared error, which is averaged
over forecast cases.  Thus, the performance criterion takes the form
\begin{equation}  \label{eq:S} 
\bar{\myS} = \frac{1}{n} \sum_{i=1}^n \myS(x_i,y_i),   
\end{equation} 
where there are $n$ forecast cases with corresponding point forecasts,
$x_1, \ldots, x_n$, and verifying observations, $y_1, \ldots, y_n$.
The function $\myS$ depends both on the forecast and the realization,
and we refer to it as a {\em scoring function}.

Table \ref{tab:scores} lists some commonly used scoring functions.  We
generally take scoring functions to be {\em negatively oriented}, that
is, the smaller, the better.  The absolute error and the squared error
are of the {\em prediction error}\/ form, in that they depend on the
forecast error, $x-y$, only, and they are {\em symmetric}, in that
$\myS(x,y) = \myS(y,x)$.  The absolute percentage error and the
relative error are used for strictly positive quantities only; they
are neither of the prediction error form nor symmetric.  Patton (2009)
discusses these as well as many other scoring functions that have been
used to assess point forecasts for a strictly positive quantity, such
as an asset value or a volatility proxy.

\begin{table}[h]

\medskip
\caption{Some commonly used scoring functions. \label{tab:scores}}

\small

\begin{center}
\begin{tabular}{lll} 
\hline
\hline
$\myS(x,y) = (x-y)^2$ \hspace{5mm} & squared error (SE) \rule{0mm}{4mm} \\ 
$\myS(x,y) = \hsp |x-y|$           & absolute error (AE) \rule{0mm}{3.5mm} \\ 
$\myS(x,y) = \hsp |(x-y)/y \hsp|$  & absolute percentage error (APE) \rule{0mm}{3.5mm} \\ 
$\myS(x,y) = \hsp |(x-y)/x \hsp|$  & relative error (RE) \rule{0mm}{3.5mm} \\ 
\hline
\end{tabular}
\end{center} 

\end{table}  

Our next two tables summarize the use of scoring functions in
academia, the public and the private sector.  Table \ref{tab:lit}
surveys the 2008 volumes of peer-reviewed journals in forecasting
(Group I) and statistics (Group II), along with premier journals in
the most prominent application areas, namely econometrics (Group III)
and meteorology (Group IV).  We call an article a {\em forecasting
paper}\/ if it contains a table or a figure in which the predictive
performance of a forecaster or forecasting method is summarized in the
form of the mean score (\ref{eq:S}), or a monotone transformation
thereof, such as the root mean squared error.  Not surprisingly, the
majority of the Group I papers are forecasting papers, and many of
them employ several scoring functions simultaneously.  Overall, the
squared error is the most popular scoring function in academia,
particularly in Groups III and IV, followed by the absolute error and
the absolute percentage error.

Table \ref{tab:business} reports the use of scoring functions in
businesses and organizations, according to surveys conducted or
summarized by Carbone and Armstrong (1982), Mentzner and Kahn (1995),
McCarthy et al.~(2006) and Fildes and Goodwin (2007).  In addition to
the squared error and the absolute error, the absolute percentage
error has been very widely used in practice, presumably because
business forecasts focus on demand, sales, or costs, all of which are
nonnegative quantities.

\begin{table}[p]

\caption{Use of scoring functions in the 2008 volumes of leading
  peer-reviewed journals in forecasting (Group I), statistics (Group
  II), econometrics (Group III) and meteorology (Group IV).  Column 2
  shows the total number of papers published in 2008 under Web of
  Science document type article, note or review.  Column 3 shows the
  number of forecasting papers (FP), that is, the number of articles
  with a table or figure that summarizes predictive performance in the
  form of the mean score (\ref{eq:S}) or a monotone transformation
  thereof.  Columns 4 through 7 show the number of papers employing
  the squared error (SE), absolute error (AE), absolute percentage
  error (APE), or miscellaneous (MSC) other scoring functions.  The
  sum of columns 4 through 7 may exceed the number in column 3,
  because of the simultaneous use of multiple scoring functions in
  some articles.  Papers that apply error measures to evaluate
  estimation methods, rather than forecasting methods, have not been
  considered in this study.
  \label{tab:lit}}

\small

\begin{center} 
\begin{tabular}{lrrrrrr} 
\hline
\hline
\rule{0mm}{4.5mm} & Total & FP & SE & AE & APE & MSC \\
\hline
Group I: {\em Forecasting} \rule{0mm}{4.5mm} &&&&&& \\
\hline
International Journal of Forecasting \rule{0mm}{3.5mm}            & 41 & 32 & 21 & 10 & 8 & 4 \\ 
Journal of Forecasting                                            & 39 & 25 & 23 & 13 & 5 & 3 \\
\hline
Group II: {\em Statistics} \rule{0mm}{4.5mm} &&&&& \\
\hline 
Annals of Applied Statistics \rule{0mm}{3.5mm}                    &  62 &  8 & 6 & 3 & 1 & 0 \\
Annals of Statistics                                              & 100 &  5 & 3 & 2 & 0 & 0 \\
Journal of the American Statistical Association                   & 129 & 10 & 9 & 1 & 0 & 0 \\
Journal of the Royal Statistical Society Ser.~B                   &  49 &  5 & 4 & 1 & 0 & 0 \\
\hline
Group III: {\em Econometrics} \rule{0mm}{4.5mm} &&&&& \\
\hline 
Journal of Business and Economic Statistics \rule{0mm}{3.5mm}     &  26 &  9 &  8 &  2 & 1 & 0 \\ 
Journal of Econometrics                                           & 118 &  5 &  5 &  0 & 0 & 0 \\
\hline
Group IV: {\em Meteorology} \rule{0mm}{4.5mm} &&&&& \\
\hline 
Bulletin of the American Meteorological Society \rule{0mm}{3.5mm} &  73 &  1 &  1 &  0 & 0 & 0 \\ 
Monthly Weather Review                                            & 300 & 63 & 58 &  8 & 2 & 0 \\
Quarterly Journal of the Royal Meteorological Society             & 148 & 19 & 19 &  0 & 0 & 0 \\
Weather and Forecasting                                           &  79 & 26 & 20 & 11 & 0 & 1 \\
\hline
\end{tabular}
\end{center}

\end{table}  

\begin{table}[t]

\caption{Use of scoring functions in the evaluation of point forecasts
  in businesses and organizations.  Columns 2 through 4 show the
  percentage of survey respondents using the squared error (SE),
  absolute error (AE) and absolute percentage error (APE), with the
  source of the survey listed in column 1.
  \label{tab:business}}

\small

\begin{center} 
\begin{tabular}{lrrrrr} 
\hline
\hline
Source \rule{0mm}{4.5mm} & SE & AE & APE \\
\hline
Carbone and Armstrong (1982), Table 1 \rule{0mm}{4.5mm}  & 27\% & 19\% &  9\% \\ 
Mentzner and Kahn (1995), Table VIII                     & 10\% & 25\% & 52\% \\
McCarthy, Davis, Golicic and Mentzner (2006), Table VIII &  6\% & 20\% & 45\% \\
Fildes and Goodwin (2007), Table 5                       &  9\% & 36\% & 44\% \\
\hline
\end{tabular}
\end{center}

\end{table}  

There are many options and considerations in choosing a scoring
function.  What scoring function ought to be used in practice?  Do the
standard choices have theoretical support?  Arguably, there is
considerable contention in the scientific community, along with a
critical need for theoretically principled guidance.  Some 20 years
ago, Murphy and Winkler (1987, p.~1330) commented on the state of the
art in forecast evaluation, noting that
\begin{quote} 
\small ``[\ldots] verification measures have tended to proliferate,
with relatively little effort being made to develop general concepts
and principles [\ldots] This state of affairs has impacted the
development of a science of forecast verification.''
\end{quote} 
Nothing much has changed since.  Armstrong (2001) called for further
research, while Moskaitis and Hansen (2006) asked
\begin{quote} 
\small ``Deterministic forecasting and verification: A busted
system?''  
\end{quote} 
Similarly, the recent review by Fildes et al.~(2008, p.~1158) states that
\begin{quote} 
\small ``Defining the basic requirements of a good error measure is
still a controversial issue.''
\end{quote} 

\begin{figure}[t] 

\begin{center}
\includegraphics[width=0.90\textwidth]{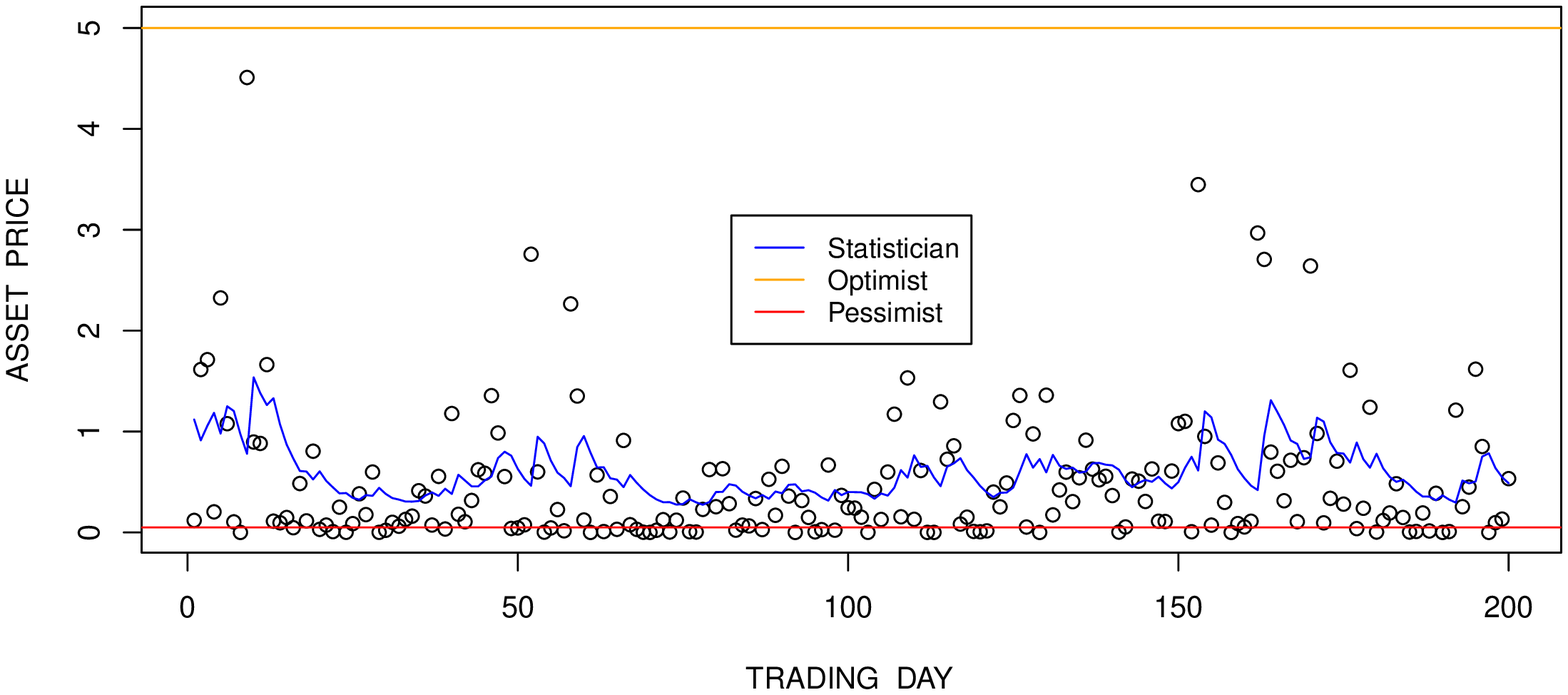}
\end{center}

\vspace{-8mm}
\caption{A realized series of volatile daily asset prices under the
  data generating process (\ref{eq:asset}), shown by circles, along
  with the one-day ahead point forecasts by the statistician (blue
  line), the optimist (orange line at top) and the pessimist (red line
  at bottom). \label{fig:sim}}

\end{figure} 

\subsection{Simulation study}  \label{sec:simulation} 

To focus issues and ideas, we consider a simulation study, in which we
seek point forecasts for a highly volatile daily asset value, $y_t$.
The data generating process is such that $y_t$ is a realization of the
random variable
\begin{equation}  \label{eq:asset} 
Y_t = Z_t^2, 
\end{equation} 
where $Z_t$ follows a Gaussian conditionally heteroscedastic time
series model (Engle 1982; Bollerslev 1986), with the parameter values
proposed by Christoffersen and Diebold (1996), in that
\[
Z_t \sim \cN(0,\sigma_t^2) 
\quad \mbox{where} \quad
\sigma_t^2 = 0.20 \, Z_{t-1}^2 + 0.75 \, \sigma_{t-1}^2 + 0.05. 
\]
We consider three forecasters, each of whom issues a one-day ahead
point forecast for the asset value.  The {\bf statistician} has
knowledge of the data generating process and the actual value of the
conditional variance $\sigma_t$, and thus predicts the true
conditional mean,
\[
\hat{x}_t = \myE \hsp (Y_t|\sigma_t^2) = \sigma_t^2,  
\]
as her point forecast.  The {\bf optimist} always predicts $\hat{x}_t
= 5$.  The {\bf pessimist} always issues the point forecast $\hat{x}_t
= 0.05$.  Figure \ref{fig:sim} shows these point forecasts along with
the realizing asset value for 200 successive trading days.  There
ought to be little contention as to the predictive performance, in
that the statistician is more skilled than the optimist or the
pessimist.

Table \ref{tab:sim1} provides a formal evaluation of the three
forecasters for a sequence of $n = 100,000$ sequential forecasts,
using the mean score (\ref{eq:S}) and the scoring functions listed in
Table \ref{tab:scores}.  The results are counterintuitive and
disconcerting, in that the pessimist has the best (lowest) score both
under the absolute error and the absolute percentage error scoring
functions.  In terms of relative error, the optimist performs best.
Yet, what we have done here is common practice in academia and
businesses, in that point forecasts are evaluated by means of these
scoring functions.

\begin{table}[t]

\caption{The mean error measure (\ref{eq:S}) for the three point
  forecasters in the simulation study, using the squared error (SE),
  absolute error (AE), absolute percentage error (APE) and relative
  error (RE) scoring functions. \label{tab:sim1}}

\small

\begin{center}
\begin{tabular}{lrcrr}  
\hline
\hline
Forecaster \rule{0mm}{4.5mm} & SE & AE & APE & RE \\
\hline
Statistician \rule{0mm}{4.5mm} &  5.07 & 0.97 &  $2.58 \times 10^5$ &  0.97 \\
Optimist \rule{0mm}{3.5mm}     & 22.73 & 4.35 & $13.96 \times 10^5$ &  0.87 \\
Pessimist \rule{0mm}{3.5mm}    &  7.61 & 0.96 &  $0.14 \times 10^5$ & 19.24 \\
\hline
\end{tabular} 
\end{center}

\end{table} 

\subsection{Discussion}  \label{sec:interpretation}

The source of these disconcerting results is aptly explained in a
recent paper by Engelberg, Manski and Williams (2009, p.~30):
\begin{quote} 
\small
``Our concern is prediction of real-valued outcomes such as firm
profit, GDP, growth, or temperature.  In these cases, the users of
point predictions sometimes presume that forecasters report the means
of their subjective probability distributions; that is, their best
point predictions under square loss.  However, forecasters are not
specifically asked to report subjective means.  Nor are they asked to
report subjective medians or modes, which are best predictors under
other loss functions.  Instead, they are simply asked to `predict' the
outcome or to provide their `best prediction', without definition of
the word `best.'  In the absence of explicit guidance, forecasters may
report different distributional features as their point predictions.
Some may report subjective means, others subjective medians or modes,
and still others, applying asymmetric loss functions, may report
various quantiles of their subjective probability distributions.''
\end{quote}

Similarly, Murphy and Daan (1985, p.~391) noted that
\begin{quote} 
\small
``It will be assumed here that the forecasters receive a `directive'
concerning the procedure to be followed [\ldots]~and that it is
desirable to choose an evaluation measure that is consistent with this
concept.  An example may help to illustrate this concept.  Consider a
continuous [\ldots]~predictand, and suppose that the directive states
`forecast the expected (or mean) value of the variable.'  In this
situation, the mean square error measure would be an appropriate
scoring rule, since it is minimized by forecasting the mean of the
(judgemental) probability distribution.  Measures that correspond with
a directive in this sense will be referred to as consistent scoring
rules (for that directive).''
\end{quote} 

Despite these well-argued perspectives, there has been little
recognition that the common practice of requesting `some' point
forecast, and then evaluating the forecasters by using `some' (set of)
scoring function(s), is not a meaningful endeavor.  In this paper, we
develop the perspectives of Murphy and Daan (1985) and Engelberg et
al.~(2009) and argue that effective point forecasting depends on
`guidance' or `directives', which can be given in one of two
complementary ways, namely, by disclosing the {\bf scoring function}
ex ante to the forecaster, or by requesting a specific {\bf
functional} of the forecaster's predictive distribution, such as the
mean or a quantile.

As to the first option, the a priori disclosure of the scoring
function allows the forecaster to tailor the point predictor to the
scoring function at hand.  In particular, this permits our
statistician forecaster to mutate into {\bf Mr.~Bayes}, who issues the
optimal point forecast, namely the Bayes rule,
\begin{equation}  \label{eq:Bayes} 
\textstyle
\hat{x} = \arg \min_x \hsp \myE_F \, \myS(x,Y),
\end{equation} 
where the expectation is taken with respect to the forecaster's
subjective or objective predictive distribution, $F$.  For example, if
the scoring function $\myS$ is the squared error, the optimal point
forecast is the mean of the predictive distribution.  In the case of
the absolute error, the Bayes rule is any median of the predictive
distribution.  The class
\begin{equation}  \label{eq:med.beta}
\myS_\beta(x,y) = \left| 1 -  \left( \frac{y}{x} \right)^\beta \right|
\qquad (\beta \not= 0) 
\end{equation}
of scoring functions nests both the absolute percentage error ($\beta
= -1$) and the relative error ($\beta = 1$) scoring functions.  If the
predictive distribution $F$ has density $f$ on the positive half-axis
and a finite fractional moment of order $\beta$, the optimal point
forecast under the loss or scoring function (\ref{eq:med.beta}) is the
median of a random variable whose density is proportional to $y^\beta
f(y)$.  We call this the {\em $\beta$-median}\/ of the probability
distribution $F$ and write ${\rm med}^{(\beta)}(F)$.  The traditional
median arises in the limit as $\beta \to 0$.

Table \ref{tab:sim2} summarizes our discussion, in that it shows the
optimal point forecast, or Bayes rule, under the scoring functions in
Table \ref{tab:scores}, both in full generality and in the special
case of the true predictive distribution under the data generating
process (\ref{eq:asset}).  Table \ref{tab:sim3} shows the mean score
(\ref{eq:S}) for the new competitor Mr.~Bayes in the simulation study,
who issues the optimal point forecast.  As expected, Mr.~Bayes
outperforms his colleagues.

\begin{table}[t] 

\caption{Bayes rules under the scoring functions in Table
  \ref{tab:scores} as a functional of the forecaster's predictive
  distribution, $F$.  The functional ${\rm med}^{(\beta)}(F)$ is
  defined in the text.  The final column specializes to the true
  predictive distribution under the data generating process
  (\ref{eq:asset}) in the simulation study.  The entry for the
  absolute percentage error (APE) is to be understood as follows.  The
  predictive distribution $F$ has infinite fractional moment of order
  $-1$, and thus ${\rm med}^{(-1)}(F)$ does not exist.  However, it is
  readily seen that the smaller the (strictly positive) point
  forecast, the smaller the expected APE.  Thus, a prudent forecaster
  will issue some very small $\epsilon > 0$ as point predictor.
  \label{tab:sim2}}

\small

\begin{center} 
\begin{tabular}{llr}  
\hline
\hline 
Scoring Function \rule{0mm}{4.5mm} & Bayes Rule & Point Forecast in Simulation Study \\
\hline
SE  \rule{0mm}{5.5mm} & $\hat{x} = {\rm mean}(F)$     & $\sigma_t^2$ \\
AE  \rule{0mm}{4.5mm} & $\hat{x} = {\rm median}(F)$   & $0.455 \, \sigma_t^2$ \\         
APE \rule{0mm}{4.5mm} & $\hat{x} = {\rm med}^{(-1)}(F)$ & $\varepsilon$ \\
RE  \rule{0mm}{4.5mm} & $\hat{x} = {\rm med}^{(1)}(F)$  & $2.366 \, \sigma_t^2$ \\         
\hline 
\end{tabular} 
\end{center} 

\bigskip
\caption{Continuation of Table \ref{tab:sim1}, showing the
  corresponding mean scores for the new competitor, Mr.~Bayes.  In the
  case of the APE, Mr.~Bayes issues the point forecast $\hat{x} =
  \epsilon = 10^{-10}$.  \label{tab:sim3}}

\vspace{-4mm}
\begin{center}
\begin{tabular}{lcccc}  
\hline
\hline
\rule{0mm}{4.5mm} & SE & AE & APE & RE \\
\hline
Mr.~Bayes \rule{0mm}{4.5mm} &  5.07 & 0.86 & 1.00 &  0.75 \\
\hline
\end{tabular} 
\end{center}

\end{table} 

An alternative to disclosing the scoring function is to request a
specific functional of the forecaster's predictive distribution, such
as the mean or a quantile, and to apply any scoring function that is
consistent with the functional, roughly in the following sense.

Let the interval $\myI$ be the potential range of the outcomes, such
as $\myI = \real$ for a real-valued quantity, or $\myI = (0,\infty)$
for a strictly positive quantity, and let the probability distribution
$F$ be concentrated on $\myI$.  Then a {\em scoring function} is any
mapping $\myS : \myI \times \myI \to [0,\infty)$.  A {\em functional}
is a potentially set-valued mapping $F \mapsto \myT(F) \subseteq
\myI$.  A scoring function $\myS$ is {\em consistent}\/ for the
functional $\myT$ if
\[
\myE_F \hsp [\myS(t,Y)] \leq \myE_F \hsp [\myS(x,Y)]
\]
for all $F$, all $t \in \myT(F)$ and all $x \in \myI$.  It is {\em
strictly consistent}\/ if it is consistent and equality of the
expectations implies that $x \in \myT(F)$.  Following Osband (1985)
and Lambert, Pennock and Shoham (2008), a functional is {\em
elicitable} if there exists a scoring function that is strictly
consistent for it.

\subsection{Plan of the paper}  \label{sec:plan} 

The remainder of the paper is organized as follows.  Section
\ref{sec:theory} develops the notions of consistency and elicitability
in a comprehensive way.  In addition to reviewing and unifying the
extant literature, we present original results on weighted scoring
functions that extend prior findings on optimal point forecasts, such
as those of Park and Stefanski (1998) and Patton (2010).  Section
\ref{sec:examples} turns to examples.  The mean functional, ratios of
expectations, quantiles and expectiles are elicitable.  Subject to
weak regularity conditions, a scoring function for a real-valued
predictand is consistent for the mean functional if and only if it is
a Bregman function, that is, of the form
\[
\myS(x,y) = \phi(y) - \phi(x) - \phi'(x) \hsp (y-x),  
\]
where $\phi$ is a convex function with subgradient $\phi'$ (Savage
1971).  More general and novel results apply to ratios of expectations
and expectiles.  A scoring function is consistent for the
$\alpha$-quantile if and only if it is generalized piecewise linear
(GPL) of order $\alpha \in (0,1)$, that is, of the form
\[
\myS(x,y) = (\one(x \geq y) -\alpha) \, ( \hsp g(x) - g(y)),
\]
where $\one(\cdot)$ denotes an indicator function and $g$ is
nondecreasing (Thomson 1979; Saerens 2000).  However, not all
functionals are elicitable.  Notably, the conditional value-at-risk
(CVaR) functional is not elicitable, despite its popularity as a risk
measure in financial applications.

The paper closes with a discussion in Section \ref{sec:discussion},
which makes a plea for change in the practice of point forecasting.  I
contend that in issuing and evaluating point forecasts, it is
essential that either the scoring function be specified ex ante, or an
elicitable target functional be named, such as an expectation or a
quantile, and scoring functions be used that are consistent for the
target functional.

\section{A decision-theoretic approach to the evaluation of point forecasts}  
         \label{sec:theory}

We now develop a theoretical framework for the evaluation of point
forecasts.  Towards this end, we review the more general, classical
decision-theoretic setting whose basic ingredients are as follows.
\begin{enumerate} 
\item[(a)] 
An {\em observation domain}, $\myO$, which comprises the potential
outcomes of a future observation.
\item[(b)] 
A class $\cF$ of probability measures on the observation domain $\myO$
(equipped with a suitable $\sigma$-algebra), which constitutes a
family of probability distributions for the future observation.
\item[(c)]  
An {\em action domain}, $\myA$, which comprises the potential actions of a
decision maker. 
\item[(d)] 
A {\em loss function}\/ $\myL : \myA \times \myO \to [0,\infty)$,
where $\myL(a,o)$ represents the monetary or societal cost when the
decision maker takes the action $a \in \myA$ and the observation $o
\in \myO$ materializes.  
\end{enumerate} 
Given a probability distribution $F \in \cF$ for the future
observation, the {\em Bayes act}\/ or {\em Bayes rule}\/ is any
decision $\hat{a} \in \myA$ such that
\begin{equation}  \label{eq:Bayes.act} 
\textstyle
\hat{a} = \arg \min_{\hsp a} \hsp \myE_F \, \myL(a,Y),  
\end{equation} 
where $Y$ is a random variable with distribution $F$.  Thus, if the
decision maker's assessment of the uncertain future is represented by
the probability measure $F$, and she wishes to minimize the expected
loss, her optimal decision is the Bayes act, $\hat{a}$.  In general,
Bayes acts need not exist nor be unique, but in most cases of
practical interest, Bayes rules exist, and frequently they are unique
(Ferguson 1967).

\subsection{Decision-theoretic setting}  \label{sec:decision.theory} 

Point forecasting falls into the general decision-theoretic setting,
if we assume that the observation domain and the action domain
coincide.  In what follows we assume, for simplicity, that this common
domain,
\[
\myD = \myO = \myA \subseteq \real^d,
\]
is a subset of the Euclidean space $\real^d$ and equipped with the
corresponding Borel $\sigma$-algebra.  Furthermore, we refer to the
loss function as a scoring function.  With these adaptations, the
basic components of our decision-theoretic framework are as follows.
\begin{enumerate} 
\item[(a)] 
A {\em prediction-observation}\/ {\rm (PO)} {\em domain}, $\cD = \myD
\times \myD$, which is the Cartesian product of the domain $\myD
\subseteq \real^d$ with itself.
\item[(b)] 
A family $\cF$ of potential probability distributions for the future
observation $Y$ that takes values in $\myD$.
\item[(c)] 
A {\em scoring function}\/ $\myS : \cD = \myD \times \myD \to
[0,\infty)$, where $\myS(x,y)$ represents the loss or pe- nalty when the
point forecast $x \in \myD$ is issued and the observation $y \in \myD$
materializes.
\end{enumerate} 
In this setting, the {\em optimal point forecast}\/ under the
probability distribution $F \in \cF$ for the future observation, $Y$,
is the Bayes act or Bayes rule (\ref{eq:Bayes.act}), which can now be
written as
\begin{equation}  \label{eq:Bayes.rule} 
\textstyle
\hat{x} = \arg \min_{\hsp x} \hsp \myE_F \, \myS(x,Y).
\end{equation} 
We will mostly work in dimension $d = 1$, in which any connected
domain $\myD$ is simply an interval, $\myI$.  The cases of prime
interest then are the real line, $\myI = \real$, and the nonnegative
or positive halfaxis, $\myI = [0,\infty)$ or $\myI = (0,\infty)$.

\begin{table}[t] 
\caption{Assumptions on a scoring function $\myS$ on a\/ {\rm PO}
         domain $\cD = \myI \times \myI$, where $\myI \subseteq \real$
         is an interval, $x \in \myI$ denotes the point forecast and
         $y \in \myI$ the realizing observation.
         \label{tab:assumptions}}

\small

\begin{center}
\begin{tabular}{ll} 
\hline
\hline
(S0) & $\myS(x,y) \geq 0$ with equality if $x = y$ \rule{0mm}{4.25mm} \\
(S1) & $\myS(x,y)$ is continuous in $x$ \rule{0mm}{4mm} \\
(S2) & The partial derivative $\myS_{(1)}(x,y)$ exists and is continuous in $x$ whenever $x \not= y$ 
       \rule{0mm}{4mm} \\ 
\hline
\end{tabular}
\end{center} 
 
\end{table}  

Table \ref{tab:assumptions} summarizes assumptions which some of our
subsequent results impose on scoring functions.  The nonnegativity
condition (S0) is standard and not restrictive.  Indeed, if $\myS_0$
is such that $\myS_0(x,y) \geq \myS_0(y,y)$ for all $x, y \in \myI$,
which is a natural assumption on a loss or scoring function, then
$\myS(x,y) = \myS_0(x,y) - \myS_0(y,y)$ satisfies (S0) and shares the
optimal point forecast (\ref{eq:Bayes.rule}), subject to integrability
conditions that are not of practical concern.  Generally, a loss
function can be multiplied by a strictly positive constant and any
function that depends on $y$ only can be added, without changing the
nature of the optimal point forecast.  Furthermore, the optimization
problem in (\ref{eq:Bayes.rule}) is posed in terms of the point
predictor, $x$.  In this light, it is natural that assumptions (S1)
and (S2) concern continuity and differentiability with respect to the
first argument, the point forecast $x$.

Efron (1991) and Patton (2010) argue that\/ {\em homogeneity}\/ or
{\em scale invariance}\/ is a desirable property of a scoring
function.  We adopt this notion and call a scoring function $\myS$ on
the PO domain $\cD = \myD \times \myD$ {\em homogeneous}\/ of {\em
order}\/ $b$ if
\[
\myS(c \hsp x, c \hsp y) = |c|^{\hsp b} \, \myS(x,y) 
\quad \mbox{for all} \quad 
x, y \in \myD \;\, \mbox{and} \;\; c \in \real
\]
which are such that $c \hsp x \in \myD$ and $c \hsp y \in \myD$.
Evidently, the underlying quest is that for equivariance in the
decision problem.  The scoring function $\myS$ on the PO domain $\cD =
\myD \times \myD$ is {\em equivariant}\/ with respect to some class
$\cH$ of injections $h : \myD \to \myD$ if
\[
{\textstyle \arg \min_x} \, \myE_F \! \left[ \myS(x,h(Y)) \right] 
= h \! \left( {\textstyle \arg \min_x} \hsp \myE_F \! \left[ \myS(x,Y) \right] \right)   
\]
for all $h \in \cH$ and all probability distributions $F$ that are
concentrated on $\myD$.  For instance, if $\myS$ is homogeneous on
$\myD = \real^d$ or $\myD = (0,\infty)^d$ then it is equivariant with
respect to the multiplicative group of the linear transformations $\{
x \mapsto c \hsp x : c > 0\}$.  If the scoring function is of the
prediction error form on $\myD = \real^d$, then it is equivariant with
respect to the translation group $\{ x \mapsto x + b : b \in \real^d
\}$.

While our decision-theoretic setting resembles and follows those of
Osband (1985) and Lambert et al.~(2008), and the subsequent
development owes much to their pioneering works, there are
distinctions in technique.  For example, Osband (1985) assumes a
bounded domain $\myD$, while Lambert et al.~(2008) consider $\myD$ to
be a finite set.  The work of Granger and Pesaran (2000a, 2000b),
which argues in favor of closer links between decision theory and
forecast evaluation, focuses on probability forecasts for a
dichotomous event.

\subsection{Consistency}  \label{sec:consistency} 

In the decision-theoretic framework, we think of the aforementioned
`distributional feature' or `directive' for the forecaster as a
statistical functional.  Formally, a {\em statistical functional}, or
simply a {\em functional}, is a potentially set-valued mapping from a
class of probability distributions, $\cF$, to a Euclidean space
(Horowitz and Manski 2006; Huber and Ronchetti 2009; Wellner 2009).
In the current context of point forecasting, we require that the
functional
\[
\myT : \cF \longrightarrow \myD, \qquad F \longmapsto \myT(F),
\]
maps into the domain $\myD \subseteq \real^d$.  Frequently, we take
$\cF$ to be the class of all probability measures on $\myD$, or the
class of the probability measures with compact support in $\myD$.

To facilitate the presentation, the following definitions and results
suppress the dependence of the scoring function $\myS$, the functional
$\myT$ and the class $\cF$ on the domain $\myD$.

\smallskip
\begin{definition}  \label{def:consistent} 
{\em The scoring function\/ $\myS$ is} consistent\/ {\em for the
functional\/ $\myT$ relative to the class\/ $\cF$ if
\begin{equation}  \label{eq:consistent} 
\myE_F \, \myS(t,Y) \leq \myE_F \, \myS(x,Y) 
\end{equation}
for all probability distributions\/ $F \in \cF$, all $t \in \myT(F)$
and all\/ $x \in \myD$.  It is} strictly consistent\/ {\em if it is
consistent and equality in {\rm (\ref{eq:consistent})} implies that $x
\in \myT(F)$.}
\end{definition}

As noted, the term consistent was coined by Murphy and Daan (1985,
p.~391), who stressed that is is critically important to define
consistency for a fixed, given functional, as opposed to a generic
notion of consistency, which was, correctly, refuted by Jolliffe
(2008).  For example, the squared error scoring function, $\myS(x,y) =
(x-y)^2$, is consistent, but not strictly consistent, for the mean
functional relative to the class of the probability measures on the
real line with finite first moment.  It is strictly consistent
relative to the class of the probability measures with finite second
moment.

In a parametric context, Lehmann (1951) and Noorbaloochi and Meeden
(1983) refer to a related property as {\em decision-theoretic
  unbiasedness}.  The following result notes that consistency is the
dual of the optimal point forecast property, just as
decision-theoretic unbiasedness is the dual of being Bayes
(Noorbaloochi and Meeden 1983).  It thus connects the problems of
finding optimal point forecasts, and of evaluating point predictions.

\smallskip
\begin{theorem}  \label{th:dual}
The scoring function\/ $\myS$ is consistent for the functional\/
$\myT$ relative to the class\/ $\cF$ if and only if, given any\/ $F
\in \cF$, any\/ $x \in \myT(F)$ is an optimal point forecast under\/
$\myS$.
\end{theorem} 

Stated differently, the class of the scoring functions that are
consistent for a certain functional is identical to the class of the
loss functions under which the functional is an optimal point
forecast.  Despite its simplicity, and the proof being immediate from
the defining properties, this duality does not appear to be widely
appreciated.

Our next result shows that the class of the consistent scoring
functions is convex, and thus suggests the existence of Choquet
representations (Phelps 1966). 

\smallskip
\begin{theorem}  \label{th:convex}
Let\/ $\lambda$ be a measure on a measurable space\/ $(\Omega,\cA)$.
Suppose that for all\/ $\omega \in \Omega$, the scoring function\/
$\myS_{\hsp \omega}$ satisfies {\rm (S0)} and is consistent for the
functional\/ $\myT$ relative to the class\/ $\cF$.  Then the scoring
function
\[
\myS(x,y) = \int \myS_{\hsp \omega}(x,y) \: \lambda(\dd \omega)
\]
is consistent for\/ $\myT$ relative to\/ $\cF$.  
\end{theorem} 

At this point, it will be useful to distinguish the notions of a
proper scoring rule (Winkler 1996; Gneiting and Raftery 2007) and a
consistent scoring function.  I believe that this distinction
is useful, even though the extant literature has failed to make it.
For example, in referring to proper scoring rules for quantile
forecasts, Cervera and Mu\~{n}oz (1996), Gneiting and Raftery (2007),
Hilden (2008) and Jose and Winkler (2009) discuss scoring functions
that are consistent for a quantile.  

Within our decision-theoretic framework, a {\em proper scoring rule}\/
is a function $\bS : \cF \times \myD \to \real$ such that
\begin{equation}  \label{eq:proper} 
\myE_F \, \bS(F,Y) \leq \myE_F \, \bS(G,Y)  
\end{equation}  
for all probability distributions\/ $F, G \in \cF$, where we assume
that the expectations are well-defined.  Note that $\bS$ is defined on
the Cartesian product of the class $\cF$ and the domain $\myD$.  The
loss or penalty $\bS(F,y)$ arises when a probabilistic forecaster
issues the predictive distribution $F$ while $y \in \myD$
materializes.  The expectation inequality (\ref{eq:proper}) then
implies that the forecaster minimizes the expected loss by following
her true beliefs.  Thus, the use of proper scoring rules encourages
sincerity and candor among probabilistic forecasters.

In contrast, a scoring function $\myS$ acts on the PO domain, $\cD =
\myD \times \myD$, that is, the Cartesian product of $\myD$ with
itself.  This is a much simpler domain than that for a scoring rule.
However, any consistent scoring function induces a proper scoring rule
in a straightforward and natural construction, as follows.

\smallskip
\begin{theorem} \label{th:proper}
Suppose that the scoring function $\myS$ is consistent for the
functional\/ $\myT$ relative to the class\/ $\cF$.  Then the function
\[
\bS : \cF \times \myD \longrightarrow [0,\infty), \qquad 
(F,y) \longmapsto \bS(F,y) = \myS(\myT(F),y),  
\]
is a proper scoring rule.  
\end{theorem} 

A more general decision-theoretic approach to the construction of
proper scoring rules is described by Dawid (2007, p.~78) and Gneiting
and Raftery (2007, p.~361).

\subsection{Elicitability}  \label{sec:elicitability} 

We turn to the notion of elicitability, which is a critically
important concept in the evaluation of point forecasts.  While the
general notion dates back to the pioneering work of Osband (1985), the
term elicitable was coined only recently by Lambert et al.~(2008).
Whenever appropriate and feasible, we suppress the dependence of the
definitions and results on the PO domain $\cD = \myD \times \myD$.

\smallskip
\begin{definition}  \label{def:elicitability} 
{\em The functional $\myT$ is} elicitable\/ {\em relative to the class $\cF$ if
there exists a scoring function $\myS$ that is strictly consistent
for $\myT$ relative to $\cF$.}
\end{definition} 

Evidently, if $\myT$ is elicitable relative to the class $\cF$, then
it is elicitable relative to any subclass $\cF_0 \subseteq \cF$.  The
following result then is a version of Osband's (1985, p.~9) {\em
revelation principle}.

\smallskip
\begin{theorem}[Osband]  \label{th:relevation} 
Suppose that the class\/ $\cF$ is concentrated on the domain\/ $\myD$,
and let\/ $g : \myD \to \myD$ be a one-to-one mapping.  Then the
following holds.
\begin{itemize} 
\item[\rm (a)] 
If\/ $\myT$ is elicitable, then\/ $\myT_g = g \circ \myT$ is elicitable.
\item[\rm (b)] 
If\/ $\myS$ is consistent for\/ $\myT$, then the scoring function\/ 
\[
\myS_g(x,y) = \myS( \hsp g^{-1}(x),y)
\]
is consistent for\/ $\myT_g$.
\item[\rm (c)] 
If\/ $\myS$ is strictly consistent for\/ $\myT$, then\/ $\myS_g$ is
strictly consistent for\/ $\myT_g$.
\end{itemize}
\end{theorem} 

The next theorem is an original result that concerns weighted scoring
functions, where the weight function depends on the realizing
observation, $y$, only.

\smallskip
\begin{theorem}  \label{th:weighted}
Let the functional\/ $\myT$ be defined on a class\/ $\cF$ of
probability distributions which admit a density, $f$, with respect to
some dominating measure on the domain\/ $\myD$.  Consider the weight
function
\[
w : \myD \to [0,\infty). 
\]
Let\/ $\cF^{\hsp (w)} \subseteq \cF$ denote the subclass of the
probability distributions in\/ $\cF$ which are such that\/ $w(y) \hsp
f(y)$ has finite integral over\/ $\myD$, and the probability measure\/ 
$F^{\hsp (w)}$ with density proportional to\/ $w(y) \hsp f(y)$ belongs
to\/ $\cF$.  Define the functional
\begin{equation}  \label{eq:T^w}
\myT^{\hsp (w)} : \hsp \cF^{\hsp (w)} \longrightarrow \myI,
\qquad 
F \longmapsto \myT^{\hsp (w)}(F) = \myT(F^{\hsp (w)}), 
\end{equation}  
on this subclass\/ $\cF^{\hsp (w)}$.  Then the following holds.
\begin{itemize} 
\item[\rm (a)] 
If\/ $\myT$ is elicitable, then\/ $\myT^{\hsp (w)}$ is elicitable.
\item[\rm (b)] 
If\/ $\myS$ is consistent for\/ $\myT$ relative to $\cF$, then 
\begin{equation}  \label{eq:S^w}
\myS^{(w)}(x,y) = w(y) \, \myS(x,y)
\end{equation}  
is consistent for\/ $\myT^{\hsp (w)}$ relative to\/ $\cF^{(w)}$.
\item[\rm (c)] 
If\/ $\myS$ is strictly consistent for\/ $\myT$ relative to $\cF$,
then\/ $\myS^{(w)}$ is strictly consistent for\/ $\myT^{\hsp (w)}$
relative to\/ $\cF^{(w)}$.
\end{itemize}
\end{theorem} 

In other words, a weighted scoring function is consistent for the
functional $\myT^{\hsp (w)}$, which acts on the predictive
distribution in a peculiar way, in that it applies the original
functional, $\myT$, to the probability measure whose density is
proportional to the product of the weight function and the original
density.

Theorem \ref{th:weighted} is a very general result with a wealth of
applications, both in forecast evaluation and in the derivation of
optimal point forecasts.  In particular, the functional (\ref{eq:T^w})
is the optimal point forecast under the weighted scoring function
(\ref{eq:S^w}), which allows us to unify and extend scattered prior
results.  For example, the scoring function $\myS_\beta$ of equation
(\ref{eq:med.beta}),
\[
\myS_\beta(x,y) = \left| 1 -  \left( \frac{y}{x} \right)^\beta \right| \! ,
\]
is of the form (\ref{eq:S^w}) with the original scoring function
$\myS(x,y) = |x^{-\beta} - y^{-\beta}|$ and the weight function $w(y)
= y^\beta$ on the positive halfaxis, $\myD = (0,\infty)$.  The scoring
function $\myS$ is consistent for the median functional.  Thus, as
noted in the introduction, the scoring function $\myS_\beta$ is
consistent for the {\em $\beta$-median}\/ functional, ${\rm
med}^{(\beta)}(F)$, that is, the median of a probability distribution
whose density is proportional to $y^\beta f(y)$, where $f$ is the
density of $F$.  If $\beta = - 1$, we recover the absolute percentage
error, $\myS_{-1}(x,y) = |(x-y)/y|$.  The case $\beta = 1$ corresponds
to the relative error, $\myS_1(x,y) = |(x-y)/x|$, which Patton (2010)
refers to as the MAE-prop function.  Table 1 of Patton (2010) shows
Monte Carlo based approximate values for optimal point forecasts under
this scoring function.  Theorem \ref{th:weighted} permits us to give
exact results; these are summarized in Table \ref{tab:Patton} and
differ notably from the approximations.

\begin{table}[t]

\caption{The optimal point forecast or Bayes rule
  (\ref{eq:Bayes.rule}) when the scoring function is relative error,
  $\myS(x,y) = |(x-y)/x|$, and the future quantity $Y$ can be
  represented as $Y = Z^2$, where $Z$ has a $t$-distribution with mean
  0, variance 1 and $\nu > 2$ degrees of freedom.  In the limiting
  case as $\nu \to \infty$, we take $Z$ to be standard normal.  If $Z$
  has variance $\sigma^2$ the entries need to be multiplied by this
  factor.  As opposed to the approximations in Table 1 of Patton
  (2010), which stem from numerical and Monte Carlo methods and are
  reproduced below, our results derive from Theorem \ref{th:weighted}
  and are exact.  For details see Appendix B.
  \label{tab:Patton}}

\small

\begin{center}
\begin{tabular}{lccccc}  
\hline 
\hline 
\rule{0mm}{4.5mm} & $\nu = 4$ & $\nu = 6$ & $\nu = 8$ & $\nu = 10$ & $\nu \to \infty$ \\ 
\hline 
Exact optimal point forecast \rule{0mm}{4.5mm}
                       & 3.4048 & 2.8216 & 2.6573 & 2.5801 & 2.3660 \\
Patton's approximation & 3.0962 & 2.7300 & 2.6067 & 2.5500 & 2.3600 \\ 
\hline
\end{tabular} 
\end{center}

\end{table} 

Another interesting case arises when the original scoring function
$\myS$ is the squared error, $\myS(x,y) = (x-y)^2$, which is
consistent for the mean or expectation functional.  If $\myT$ is the
mean functional, the functional $\myT^{(w)}$ of equation
(\ref{eq:T^w}) becomes
\begin{equation}  \label{eq:T^w.mean}
\myT^{(w)}(F) = \myT \hsp (F^{(w)}) = \myE_{F^{(w)}} [Z \hsp ]  
= \frac{\myE_F \hsp [\hsp Y \hsp w(Y)]}{\myE_F \hsp [\hsp w(Y)]}.   
\end{equation}  
Park and Stefanski (1998) studied optimal point forecasts in the
special case in which $\myD = (0,\infty)$ is the positive half-axis
and $w(y) = 1/y^2$, so that $\myS^{(w)}(x,y) = (x-y)^2 \hsp / \hsp
y^2$ is the squared percentage error.  By equation
(\ref{eq:T^w.mean}), the scoring function $\myS^{(w)}$ is consistent
for the functional $\myT^{(w)}(F) = \myE_F \hsp [Y^{-1}] \, / \,
\myE_F \hsp [Y^{-2}]$.  By Theorem \ref{th:dual}, this latter quantity
is the optimal point forecast under the squared percentage error
scoring function, which is the result derived by Park and Stefanski
(1998).

Situations in which the weight function depends on the point forecast,
$x$, need to be handled on a case-by-case basis.  For example, a
routine calculation shows that the squared relative error scoring
function, $\myS(x,y) = (x-y)^2 \hsp / x^2$, is consistent for the
functional
\begin{equation}  \label{eq:T.SRE} 
\myT(F) = \frac{\myE_F \hsp [Y^2]}{\myE_F \hsp [Y]}.  
\end{equation} 
Incidentally, by a special case of (\ref{eq:T^w.mean}) the
observation-weighted scoring function $\myS(x,y) = y \hsp (x-y)^2$ is
also consistent for the functional (\ref{eq:T.SRE}).  Later on in
equation (\ref{eq:Bregman.special}) we characterize the class of the
scoring functions that are consistent for this functional.

While Theorems \ref{th:relevation} and \ref{th:weighted} suggest that
general classes of functionals are elicitable, not all functionals are
such.  The following result, which is a variant of Proposition 2.5 of
Osband (1985) and Lemma 1 of Lambert et al.~(2008), states a necessary
condition.

\smallskip
\begin{theorem}[Osband]  \label{th:necessary} 
If a functional is elicitable then its level sets are convex in the
following sense: If\/ $F_0 \in \cF$, $F_1 \in \cF$ and\/ $p \in (0,1)$
are such that\/ $F_p = (1-p) F_0 + p F_1 \in \cF$, then\/ $t \in
\myT(F_0)$ and\/ $t \in \myT(F_1)$ imply\/ $t \in \myT(F_p)$.
\end{theorem} 

For example, the sum of two distinct quantiles generally does not have
convex level sets and thus is not an elicitable functional.
Interesting open questions include those for a converse of Theorem
\ref{th:necessary} and, more generally, for a characterization of
elicitability.

\subsection{Osband's principle}  \label{sec:Osband} 

Given an elicitable functional $\myT$, is there a practical way of
describing and characterizing the class of the scoring functions that
are consistent for it?  The following general approach, which
originates in the pioneering work of Osband (1985), is frequently
useful.

Suppose that the functional $\myT$ is defined for a class of
probability measures on the domain $\myD$ which includes the two-point
distributions.  Assume that there exists an {\em identification
function}\/ $\myV : \myD \times \myD \to \real$ such that
\begin{equation}  \label{eq:V} 
\myE_F \hsp [ \hsp \myV(x,Y)] = 0 \iff x \in \myT(F)
\end{equation} 
and $\myV(x,y) \not= 0$ unless $x = y$.  If a consistent scoring
function is available, which is smooth in its first argument, we can
take $\myV(x,y)$ to be the corresponding partial derivative.  For
example, if $\myT$ is the mean or expectation functional on an
interval $\myD = \myI \subseteq \real$, we can pick $\myV(x,y) = x -
y$, which derives from the squared error scoring function, $\myS(x,y)
= (x-y)^2$.  Table \ref{tab:V} provides further examples, with the
second and fourth nesting the first.

\begin{table}[t]

\caption{Possible choices for the identification function $\myV$ with the
         property (\ref{eq:V}) in the case in which $\myD = \myI
         \subseteq \real$ is an interval. \label{tab:V}}

\small

\begin{center}
\begin{tabular}{ll}  
\hline 
\hline 
Functional & Identification function \rule{0mm}{4.5mm} \\
\hline 
Mean \rule{0mm}{4.5mm} & $\myV(x,y) = x - y$ \\
Ratio $\myE_F \hsp [r(Y)] \, / \, \myE_F \hsp [s(Y)]$ & $\myV(x,y) = x \hsp s(y) - r(y)$ \\
$\alpha$-Quantile & $\myV(x,y) = \one(x \geq y) - \alpha$ \\
$\tau$-Expectile  & $\myV(x,y) = 2 \, |\one(x \geq y) - \tau \hsp | \, (x - y)$ \\
\hline
\end{tabular} 
\end{center}

\end{table} 

The function
\begin{equation}  \label{eq:Osband1} 
\epsilon(c) = p \, \myS(c,a) + (1-p) \hsp \myS(c,b)  
\end{equation} 
represents the expected score when we issue the point forecast $c$ for
a random vector $Y$ such that $Y = a$ with probability $p$ and $Y =
b$ with probability $1-p$.  Since $\myS$ is consistent for the
functional $\myT$, the identification function property (\ref{eq:V})
implies that $\epsilon(c)$ has a minimum at $c = x$, where
\begin{equation}  \label{eq:Osband2} 
p \, \myV(x,a) + (1-p) \, \myV(x,b) = 0.
\end{equation} 
If $\myS$ is smooth in its first argument, we can combine (\ref{eq:Osband1}) 
and (\ref{eq:Osband2}) to result in   
\begin{equation}  \label{eq:Osband3} 
\myS_{(1)}(x,a) \hsp / \, \myV(x,a) = \myS_{(1)}(x,b) \hsp / \, \myV(x,b),
\end{equation} 
where $\myS_{(1)}$ denotes a partial derivative or gradient with
respect to the first argument.  If this latter equality holds for all
pairwise distinct $a, b$ and $x \in \myD$, the function
$\myS_{(1)}(x,y) \hsp / \hsp \myV(x,y)$ is independent of $y \in
\myD$, and we can write
\begin{equation}  \label{eq:ansatz} 
\myS_{(1)}(x,y) = h(x) \, \myV(x,y) 
\end{equation} 
for $x, y \in \myD$ and some function $h : \myD \to \myD$.
Frequently, we can integrate (\ref{eq:ansatz}) to obtain the general
form of a scoring rule that is consistent for the functional $\myT$.

In recognition of Osband's (1985) fundamental yet unpublished work, we
refer to this general approach as {\em Osband's principle}.  The
examples in the subsequent section give various instances in which the
principle can be successfully put to work.  For a general technical
result, see Theorem 2.1 of Osband (1985).

\section{Examples}  \label{sec:examples} 

We now give examples in the case of a univariate predictand, in which
any connected domain $\myD = \myI \hsp \subseteq \real$ is an
interval.  Some of the results are classical, such as the
characterizations for expectations (Savage 1971) and quantiles
(Thomson 1979), and some are novel, including those for ratios of
expectations, expectiles and conditional value-at-risk.  In a majority
of the examples, the technical arguments rely on the properties of
convex functions and subgradients, for which we refer to Rockafellar
(1970).

\subsection{Expectations}  \label{sec:mean} 

It is well known that the squared error scoring function, $\myS(x,y) =
(x-y)^2$, is strictly consistent for the mean functional relative to
the class of the probability distributions on $\real$ whose second
moment is finite.  Thus, means or expectations are elicitable.  Before
turning to more general settings in subsequent sections, we review a
classical result of Savage (1971) which identifies the class of the
scoring functions that are consistent for the mean functional as that
of the Bregman functions.  Closely related results have been obtained
by Reichelstein and Osband (1984), Saerens (2000), Banerjee, Guo and
Wang (2005) and Patton (2010).

\smallskip
\begin{theorem}[Savage]  \label{th:mean} 
Let\/ $\cF$ be the class of the probability measures on the interval\/
$\myI \subseteq \real$ with finite first moment.  Then the following
holds.
\begin{enumerate} 
\item[\rm (a)] 
The mean functional is elicitable relative to the class\/ $\cF$. 
\item[\rm (b)] 
Suppose that the scoring function\/ $\myS$ satisfies assumptions\/
{\rm (S0)}, {\rm (S1)} and\/ {\rm (S2)} on the\/ {\rm PO} domain\/
$\cD = \myI \times \myI$.  Then\/ $\myS$ is consistent for the mean
functional relative to the class of the compactly supported
probability measures on\/ $\myI$ if, and only if, it is of the form
\begin{equation}  \label{eq:Bregman} 
\myS(x,y) = \phi(y) - \phi(x) - \phi'(x) \hsp (y-x),  
\end{equation} 
where $\phi$ is a convex function with subgradient\/ $\phi'$ on\/
$\myI$.  
\item[\rm (c)] 
If\/ $\phi$ is strictly convex, the scoring function\/ {\rm
(\ref{eq:Bregman})} is strictly consistent for the mean functional
relative to the class of the probability measures\/ $F$ on\/ $\myI$
for which both\/ $\myE_F \hsp Y$ and\/ $\myE_F \, \phi(Y)$ exist and
are finite.
\end{enumerate} 
\end{theorem} 

Banerjee et al.~(2005) refer to a function of the form
(\ref{eq:Bregman}) as a {\em Bregman function}.  For example, if $\myI
= \real$ and $\phi(x) = |x|^a$, where $a > 1$ to ensure strict
convexity, the Bregman representation yields the scoring function
\begin{equation}  \label{eq:spherical} 
\myS_a(x,y) = |y|^a - |x|^a - a \, \mbox{sign}(x) \hsp |x|^{a-1} (y - x),  
\end{equation}  
which is homogeneous of order $a$ and nests the squared error that
arises when $a = 2$.  Savage (1971) showed that up to a multiplicative
constant squared error is the unique Bregman function of the
prediction error form, as well as the unique symmetric Bregman
function.  Patton (2010) introduced a rich and flexible family of
homogeneous Bregman functions on the PO domain $\cD = (0,\infty)
\times (0,\infty)$, namely
\begin{equation}  \label{eq:Patton} 
\myS_b(x,y) = \left\{ \begin{array}{lcl} 
  \displaystyle \frac{1}{b \hsp (b-1)} \left( y^b - x^b \right) -
  \frac{1}{b-1} \, x^{b-1} \, (y-x) & \mbox{if} & b \in \real
  \setminus \{0,1\}, \\ \displaystyle \frac{y}{x} - \log\frac{y}{x} -
  1 & \mbox{if} & b = 0, \rule{0mm}{6.75mm} \\ \displaystyle y
  \log\frac{y}{x} - y + x & \mbox{if} & b = 1. \rule{0mm}{6mm} \\
\end{array} \right.
\end{equation}
Up to a multiplicative constant, these are the only homogeneous
Bregman functions on this PO domain.  The squared error scoring
function emerges when $b = 2$ and the QLIKE function (Patton 2010)
when $b = 0$.  If $b = a > 1$ the Patton function (\ref{eq:Patton})
coincides with the corresponding restriction of the power function
(\ref{eq:spherical}), up to a multiplicative constant.  

Finally, it is worth noting that roper scoring rules for probability
forecasts of a dichotomous event are also of the Bregman form, because
the probability of a binary event equals the expectation of the
corresponding indicator variable.  Compare McCarthy (1956), Savage
(1971), DeGroot and Fienberg (1983), Schervish (1989), Winkler (1996),
Buja, Stuetzle and Shen (2005) and Gneiting and Raftery (2007), among
others.

\begin{figure}[t] 

\begin{center}
\includegraphics[width=0.70\textwidth]{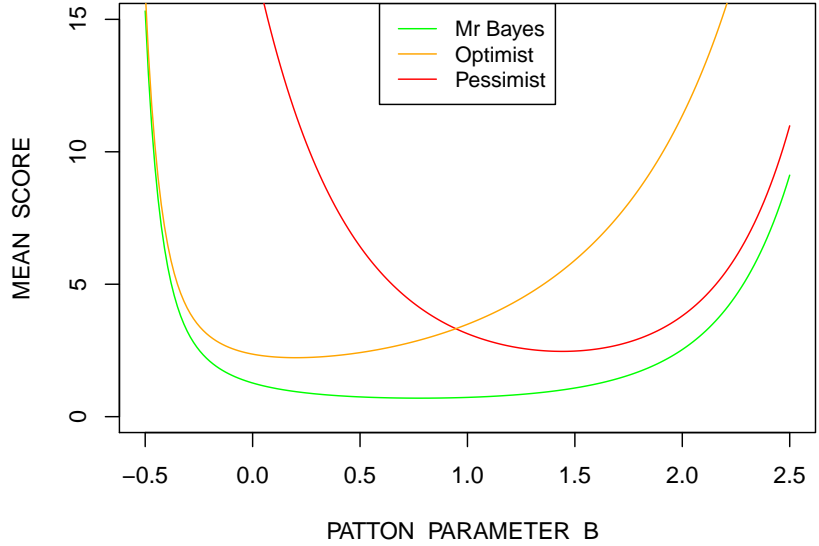}
\end{center}

\vspace{-8mm}
\caption{The mean score (\ref{eq:S}) under the Patton scoring function
  (\ref{eq:Patton}) for Mr.~Bayes (green), the optimist (orange) and
  the pessimist (red) in the simulation study of Section
  \ref{sec:simulation}.  \label{fig:Patton}}

\end{figure} 

Figure \ref{fig:Patton} returns to the initial simulation study of
Section \ref{sec:simulation} and shows the mean score (\ref{eq:S})
under the Patton scoring function (\ref{eq:Patton}) for Mr.~Bayes, the
optimist and the pessimist.  The optimal point forecast under a
Bregman scoring function is the mean of the predictive distribution,
so that the statistician forecaster fuses with Mr.~Bayes.

\subsection{Ratios of expectations}  \label{sec:ratio} 

We now consider statistical functionals which can be represented as
ratios of expectations.  The mean functional emerges in the special
case in which $r(y) = y$ and $s(y) = 1$.

\smallskip
\begin{theorem}  \label{th:ratio} 
Let $\myI \subseteq \real$ be an interval, and suppose that\/ $r :
\myI \to \real$ and\/ $s : \myI \to (0,\infty)$ are measurable
functions.  Then the following holds.
\begin{enumerate} 
\item[\rm (a)] 
The functional 
\begin{equation}  \label{eq:ratio} 
\myT(F) = \Ers\end{equation} is elicitable relative to the class of
the probability measures on\/ $\myI$ for which\/ $\myE_F \hsp [ \hsp
r(Y)]$, $\myE_F \hsp [ \hsp s(Y)]$ and\/ $\myE_F \hsp [ \hsp\hsp Y
s(Y) \hsp ]$ exist and are finite.
\item[\rm (b)] 
If\/ $\myS$ is of the form
\begin{equation}  \label{eq:Bregman.ratio} 
\myS(x,y) = s(y) \left( \phi(y) - \phi(x) \right)  
            - \phi'(x) \hsp (r(y) - x \hsp s(y)) + \phi'(y) \hsp (r(y) - y \hsp s(y)),  
\end{equation} 
where $\phi$ is a convex function with subgradient\/ $\phi'$, then it
is consistent for the functional\/ {\rm (\ref{eq:ratio})} relative to
the class of the probability measures $F$ on $\myI$ for which\/
$\myE_F \hsp [ \hsp r(Y)]$, $\myE_F \hsp [ \hsp s(Y)]$, $\myE_F \hsp [
\hsp r(Y) \hsp \phi'(Y)]$, $\myE_F \hsp [ \hsp s(Y) \hsp \phi(Y)]$
and\/ $\myE_F \hsp [ \hsp\hsp Y s(Y) \hsp \phi'(Y)]$ exist and are
finite.  If\/ $\phi$ is strictly convex, then\/ $\myS$ is strictly consistent. 
\item[(c)] 
Suppose that the scoring function\/ $\myS$ satisfies assumptions\/
{\rm (S0)}, {\rm (S1)} and\/ {\rm (S2)} on the\/ {\rm PO} domain\/
$\cD = \myI \times \myI$.  If\/ $s$ is continuous and\/ $r(y) = y \hsp
s(y)$ for\/ $y \in \myI$, then $\myS$ is consistent for the
functional\/ {\rm (\ref{eq:ratio})} relative to the class of the
compactly supported probability measures on\/ $\myI$ if, and only if,
it is of the form {\rm (\ref{eq:Bregman.ratio})}, where $\phi$ is a
convex function with subgradient\/ $\phi'$.
\end{enumerate} 
\end{theorem} 

In the case in which $s(y) = w(y)$ and $r(y) = y \hsp w(y)$ for a
strictly positive, continuous weight function $w$, the ratio
(\ref{eq:ratio}) coincides with the functional (\ref{eq:T^w.mean}).
If $\myI = (0,\infty)$ and $w(y) = y$, the special case $\myT(F) =
\myE_F \hsp [Y^2] \, / \, \myE_F \hsp [Y]$ of equation
(\ref{eq:T.SRE}) arises.  In Section \ref{sec:elicitability} we saw
that both the squared relative error scoring function, $\myS(x,y) =
(x-y)^2 \hsp / \hsp x^2$, and the observation-weighted scoring
function $\myS(x,y) = y \hsp (x-y)^2$ are consistent for this
functional.  By part (c) of Theorem \ref{th:ratio}, the general form
of a scoring function that is consistent for the functional
(\ref{eq:T.SRE}) is
\begin{equation}  \label{eq:Bregman.special} 
\myS(x,y) = y \left( \phi(y) - \phi(x) \right) - y \, (y-x) \, \phi'(x),  
\end{equation} 
where $\phi$ is convex with subgradient $\phi'$.  The above scoring
functions emerge when $\phi(y) = 1/y$ and $\phi(y) = y^2$,
respectively.

\subsection{Quantiles and expectiles}  \label{sec:quantiles} 

An $\alpha$-quantile $(0 < \alpha < 1)$ of the cumulative distribution
function $F$ is any number $x$ for which $\lim_{y \uparrow x} F(y)
\leq \alpha \leq F(x)$.  In finance, quantiles are often referred to
as {\em value-at-risk}\/ (VaR; Duffie and Pan 1997).  The literature
on the evaluation of quantile forecasts generally recommends the use
of the {\em asymmetric piecewise linear}\/ scoring function,
\begin{equation}  \label{eq:PL}
\myS_\alpha(x,y) = \left( \one(x \geq y) - \alpha \right) (x-y), 
\end{equation} 
which is strictly consistent for the $\alpha$-quantile relative to the
class of the probability measures with finite first moment (Raiffa and
Schlaifer 1961, p.~196; Ferguson 1967, p.~51).  This well-known
property lies at the heart of quantile regression (Koenker and Bassett
1978).

As regards the characterization of the scoring functions that are
consistent for a quantile, results of Thomson (1979) and Saerens
(2000) can be summarized as follows.  For a discussion of their
equivalence and historical comments, see Gneiting (2010).

\smallskip
\begin{theorem}[Thomson, Saerens]  \label{th:quantiles} 
Let\/ $\cF$ be the class of the probability measures on the interval\/
$\myI \subseteq \real$, and let\/ $\alpha \in (0,1)$.  Then the
following holds.
\begin{enumerate} 
\item[\rm (a)] 
The $\alpha$-quantile functional is elicitable relative to the class\/ $\cF$.
\item[\rm (b)] 
Suppose that the scoring function\/ $\myS$ satisfies assumptions\/
{\rm (S0)}, {\rm (S1)} and {\rm (S2)} on the\/ {\rm PO} domain\/ $\cD
= \myI \times \myI$.  Then\/ $\myS$ is consistent for the
$\alpha$-quantile relative to the class of the compactly supported
probability measures on\/ $\myI$ if, and only if, it is of the form
\begin{equation}  \label{eq:GPL} 
\myS(x,y) = (\one(x \geq y) - \alpha) \, ( \, g(x) - g(y)), 
\end{equation}
where\/ $g$ is a nondecreasing function on\/ $\myI$.  
\item[\rm (c)] 
If\/ $g$ is strictly increasing, the scoring function\/ {\rm
(\ref{eq:GPL})} is strictly consistent for the\/ $\alpha$-quantile 
relative to the class of the probability measures\/ $F$ on\/ $\myI$
for which\/ $\myE_F \hsp g(Y)$ exists and is finite.
\end{enumerate} 
\end{theorem} 

Gneiting (2008b) refers to a function of the form (\ref{eq:GPL}) as
{\em generalized piecewise linear}\/ (GPL) of order\/ $\alpha \in
(0,1)$, because it is piecewise linear after applying a nondecreasing
transformation.  Any GPL function is equivariant with respect to the
class of the nondecreasing transformations, just as the quantile
functional is equivariant under monotone mappings (Koenker 2005,
p.~39).  If $\myI = (0,\infty)$ and $g(x) = x^b/|b|$ for $b \in \real
\setminus \{ 0 \}$, and taking the corresponding limit as $b \to 0$,
we obtain the family
\begin{equation}  \label{eq:GPL.b} 
\myS_{\alpha, b}(x,y) = \left\{ \begin{array}{lcl} 
\displaystyle 
\left( \one(x \geq y) - \alpha \right) \hsp \frac{1}{|b \hsp |} \left( x^b - y^b \right)
& \mbox{if} & b \in \real \setminus \{0\}, \\ 
\displaystyle \left( \one(x \geq y) - \alpha \right) \hsp \log \frac{x}{y}
& \mbox{if} & b = 0, \rule{0mm}{6.75mm} 
\end{array} \right.
\end{equation}
of the GPL power scoring functions, which are homogeneous of order
$b$.  The asymmetric piecewise linear function (\ref{eq:PL}) arises
when $b = 1$, and the MAE-LOG and MAE-SD functions described by Patton
(2009) emerge when $\alpha = \frac{1}{2}$, and $b = 0$ and $b =
\frac{1}{2}$, respectively.

\begin{figure}[t] 

\begin{center}
\includegraphics[width=0.70\textwidth]{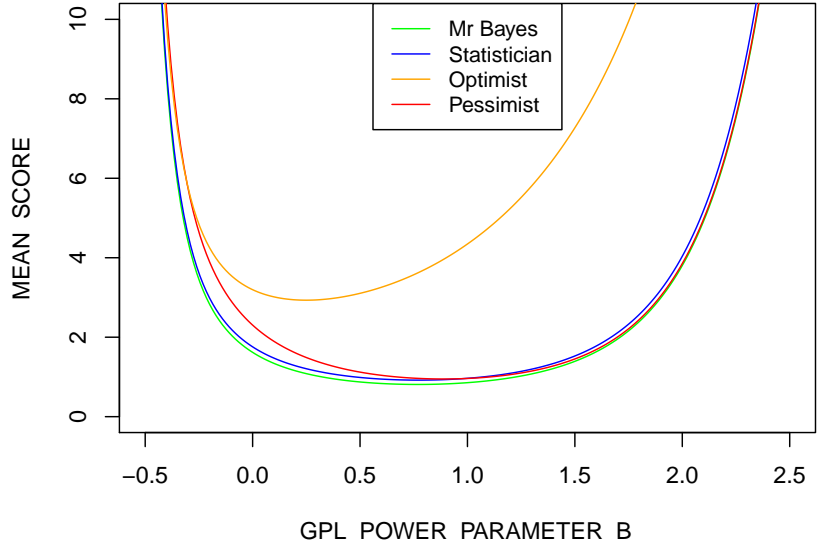}
\end{center}

\vspace{-8mm}
\caption{The mean score (\ref{eq:S}) under the GPL power scoring
  function (\ref{eq:GPL.b}) with $\alpha = \frac{1}{2}$ for Mr.~Bayes
  (green), the statistician (blue), the optimist (orange) and the
  pessimist (red) in the simulation study of Section
  \ref{sec:simulation}.  \label{fig:GPL}}

\end{figure} 

Figure \ref{fig:GPL} returns to the simulation study in Section
\ref{sec:simulation} and shows the mean score (\ref{eq:S}) under the
GPL power function (\ref{eq:GPL.b}), where $\alpha = \frac{1}{2}$, for
Mr.~Bayes, the statistician, the optimist and the pessimist.  Once
again, Mr.~Bayes dominates his competitors.

Newey and Powell (1987) introduced the {\em $\tau$-expectile}\/
functional ($0 < \tau < 1$) of a probability measure $F$ with finite
mean as the unique solution $x = \mu_\tau$ to the equation
\[
\tau \int_x^\infty (y-x) \: \dd F(y) = (1 - \tau) \int_{-\infty}^{x} (x-y) \: \dd F(y).
\]
If the second moment of $F$ is finite, the $\tau$-expectile equals the
Bayes rule or optimal point forecast (\ref{eq:Bayes.rule}) under the
{\em asymmetric piecewise quadratic}\/ scoring function,
\begin{equation}  \label{eq:tau.SE} 
\myS_\tau(x,y) = |\one(x \geq y) - \tau| \, (x-y)^2, 
\end{equation} 
similarly to the $\alpha$-quantile being the Bayes rule under the
asymmetric piecewise linear function (\ref{eq:PL}).  Not surprisingly,
expectiles have properties that resemble those of quantiles.

The following original result characterizes the class of the scoring
functions that are consistent for expectiles.  It is interesting to
observe the ways in which the corresponding class
(\ref{eq:Bregman.tau}) combines key characteristics of the Bregman and
GPL families.

\smallskip
\begin{theorem}  \label{th:expectiles} 
Let\/ $\cF$ be the class of the probability measures on the interval\/
$\myI \subseteq \real$ with finite first moment, and let\/ $\tau \in
(0,1)$.  Then the following holds.
\begin{enumerate} 
\item[\rm (a)] 
The $\tau$-expectile functional is elicitable relative to the class\/ $\cF$.
\item[\rm (b)] 
Suppose that the scoring function\/ $\myS$ satisfies assumptions\/
{\rm (S0)}, {\rm (S1)} and\/ {\rm (S2)} on the\/ {\rm PO} domain\/
$\cD = \myI \times \myI$.  Then\/ $\myS$ is consistent for the
$\tau$-expectile relative to the class of the compactly supported
probability measures on\/ $\myI$ if, and only if, it is of the form
\begin{equation}  \label{eq:Bregman.tau} 
\myS(x,y) =
|\one( x \geq y) - \tau| \left( \phi(y) - \phi(x) - \phi'(x) \hsp (y-x) \right) ,
\end{equation} 
where $\phi$ is a convex function with subgradient\/ $\phi'$ on\/
$\myI$.  
\item[\rm (c)] 
If\/ $\phi$ is strictly convex, the scoring function\/ {\rm
(\ref{eq:Bregman.tau})} is strictly consistent for the\/
$\tau$-expectile relative to the class of the probability measures\/
$F$ on\/ $\myI$ for which both\/ $\myE_F \hsp Y$ and\/ $\myE_F \hsp
\phi(Y)$ exist and are finite.
\end{enumerate} 
\end{theorem}

\subsection{Conditional value-at-risk}  \label{sec:CVaR} 

The {\em $\alpha$-conditional value-at-risk}\/ functional
($\mbox{CVaR}_{\hsp \alpha}$, $0 < \alpha < 1$) equals the expectation
of a random variable with distribution $F$ conditional on it taking
values in its upper $(1-\alpha)$-tail (Rockafellar and Uryasev 2000,
2002).  An often convenient, equivalent definition is
\begin{equation}  \label{eq:CVaR} 
{\rm CVaR}_\alpha(F) = \frac{1}{1-\alpha} \int_\alpha^1 q_\beta(F) \, \dd \beta, 
\end{equation} 
where $q_\beta$ denotes the $\beta$-quantile (Acerbi 2002), similarly
to the functional representation of the $\alpha$-trimmed mean (Huber
and Ronchetti 2009).  The CVaR functional is a popular risk measure in
quantitative finance.  Its varied, elegant and appealing properties
include coherency in the sense of Artzner et al.~(1999), who consider
functionals defined in terms of random variables, rather than the
corresponding probability measures.

\smallskip 
\begin{theorem}  \label{th:CVaR} 
The\/ $\mbox{\rm CVaR}_\alpha$ functional is not elicitable relative
to any class\/ $\cF$ of probability distributions on the interval\/
$\myI \subseteq \real$ that contains the measures with finite support,
or the finite mixtures of the absolutely continuous distributions with
compact support.
\end{theorem} 

This negative result challenges the use of the CVaR functional as a
predictive measure of risk, and may provide a partial explanation for
the striking lack of literature on the evaluation of CVaR forecasts,
as opposed to quantile or VaR forecasts, for which we refer to
Berkowitz and O'Brien (2002), Giacomini and Komunjer (2005) and Bao,
Lee and Salto\u{g}lu (2006), among others.  With consistent scoring
functions not being available, it remains unclear how one might assess
and compare CVaR forecasts.

\subsection{Mode}  \label{sec:mode} 

Let $\cF$ be a class of probability measures on the real line, each of
which has a well-defined, unique mode.  It is sometimes stated
informally that the mode is an optimal point forecast under the
zero-one scoring function,
\[
\myS_c(x,y) =  \one(|x-y| > c), 
\]
where $c > 0$.  A rigorous statement is that the optimal point
forecast or Bayes rule (\ref{eq:Bayes.rule}) under the scoring
function $\myS_c$ is the {\em midpoint}
\[
\hat{x} = {\textstyle \arg\max_{\, x}} \left( F(x+c) 
          - {\textstyle \lim_{y \uparrow x - c}} \hsp F(y) \right) 
\]
of the {\em modal interval}\/ of length $2c$ of the probability
measure $F \in \cF$ (Ferguson 1967, p.~51).  Example 7.20 of Lehmann
and Casella (1998) explores this argument in more detail.

Expressed differently, the zero-one scoring function $\myS_c$ is
consistent for the midpoint functional, which we denote by $\myT_{\!
\hsp c}$.  If $c$ is sufficiently small, then $\myT_{\! \hsp c}(F)$ is
well-defined and single-valued for all $F \in \cF$.  We can then
define the {\em mode functional}\/ on $\cF$ as the limit
\[
\myT_0(F) = {\textstyle \lim_{c \downarrow 0}} \hsp \myT_{\! \hsp c}(F). 
\]  
I do not know whether or not $\myT_0$ is elicitable.  However, if the
members of the class $\cF$ have continuous Lebesgue densities, then
$\myT_0$ is {\em asymptotically elicitable}, in the sense that it can
be represented as the continuous limit of a family of elicitable
functionals.

Stronger results become available if one puts conditions on both the
scoring function $\myS$ and the family $\cF$ of probability
distributions.  Theorem 2 of Granger (1969) is a result of this type.
Consider the PO domain $\cD = \real \times \real$.  If the scoring
function $\myS$ is an even function of the prediction error that
attains a minimum at the origin, and each $F \in \cF$ admits a
Lebesgue density, $f$, which is symmetric, continuous and unimodal, so
that mean, median and mode coincide, then $\myS$ is consistent for
this common functional.  Theorem 1 of Granger (1969) and Theorem 7.15
of Lehmann and Casella (1998) trade the continuity and unimodality
conditions on $f$ for an additional assumption of convexity on the
scoring function.

Henderson, Jones and Stare (2001, p.~3087) posit that in survival
analysis a loss function of the form
\[
\myS_k^*(x,y) = \left\{ \begin{array}{ll} 
                0, & \displaystyle \frac{x}{k} \leq y \leq kx \\
                1, & {\rm otherwise} \rule{0mm}{5mm} \end{array} \right\}  
              = \one(|\log(x) - \log(y)| > \log(k)) 
\]
is reasonable, with a choice of $k = 2$ often being adequate, arguing
that ``most people for example would accept that a lifetime prediction
of, say, 2 months, was reasonably accurate if death occurs between
about 1 and 4 months''.  From the above, the optimal point forecast or
Bayes rule under $\myS_k^*$ is the midpoint functional
$\myT_{\log(k)}$ applied to the predictive distribution of the
logarithm of the lifetime, rather than the lifetime itself.  Henderson
et al.~(2001) give various examples.

\section{Multivariate predictands} \label{sec:multivariate} 

While thus far we have restricted attention to point forecasts of a
univariate quantity, the general case of a multivariate predictand
that takes values in a domain $\myD \subseteq \real^d$ is of
considerable interest.  Applications include those of Gneiting et
al.~(2008) and Hering and Genton (2010) to predictions of wind
vectors, or that of Laurent, Rombouts and Violante (2009) to forecasts
of multivariate volatility, to name but a few.  We turn to the
decision-theoretic setting of Section \ref{sec:decision.theory} and
assume, for simplicity, that the point forecast, the observation and
the target functional take values in $\myD = \real^d$.

We first discuss the mean functional.  Assuming that $\myS(x,y) \geq
0$ with equality if $x = y$, Savage (1971), Osband and Reichelstein
(1985) and Banerjee et al.~(2005) showed that a scoring function under
which the (component-wise) expectation of the predictive distribution
is an optimal point forecast, is of the Bregman form
\begin{equation}  \label{eq:Bregman.vector} 
\myS(x,y) = \phi(y) - \phi(x) - \langle \nabla \phi(x), y-x \rangle,  
\end{equation} 
where $\phi : \real^d \to \real$ is convex with gradient $\nabla \phi
: \real^d \to \real^d$ and $\langle \: , \, \rangle$ denotes a scalar
product, subject to smoothness conditions.  Expressed differently, a
sufficiently smooth scoring function is consistent for the mean
functional if and only if it is of the form (\ref{eq:Bregman.vector}),
which is a generalization of the Bregman representation
(\ref{eq:Bregman}) in the case of a univariate predictand.  When
$\phi(x) = \|x\|^2$ is the squared Euclidean norm, we obtain the
squared error scoring function, and similarly its ramifications, such
as the weighted squared error and the pseudo Mahalanobis error
(Laurent et al.~2009).

It is of interest to note that rigorous versions of the Bregman
characterization depend on restrictive smoothness conditions.  Osband
and Reichelstein (1985) assume that the scoring function is
continuously differentiable with respect to its first argument, the
point forecast; Banerjee et al.~(2005) assume the existence of
continuous second partial derivatives with respect to the observation.
A challenging, nontrivial problem is to unify and strengthen these
results, both in univariate and multivariate settings.

Laurent et al.~(2009) consider point forecasts of multivariate
stochastic volatility, where the predictand is a symmetric and
positive definite matrix in $\real^{q \times q}$.  If the matrix is
vectorized, the above results for the mean functional apply, thereby
leading to the Bregman representation (\ref{eq:Bregman.vector}) for
the respective consistent scoring functions, which is hidden in
Proposition 3 of Laurent et al.~(2009).  Corollary 1 of Laurent et
al.~(2009) supplies a version thereof that applies directly to point
forecasts, say $\Sigma_x \in \real^{q \times q}$, of a matrix-valued,
symmetric and positive definite quantity, say $\Sigma_y \in \real^{q
\times q}$, without any need to resort to vectorization.
Specifically, any scoring function of the form
\begin{equation}  \label{eq:Bregman.matrix} 
\myS(\Sigma_x,\Sigma_y) = \phi(\Sigma_y) - \phi(\Sigma_x) - {\rm
  tr} \left( \nabla_0 \, \phi \hsp (\Sigma_x) \hsp (\Sigma_y-\Sigma_x) \right)
\end{equation} 
is consistent for the (component-wise) mean functional, where $\phi$
is convex and smooth, and $\nabla_0 \, \phi$ denotes a symmetric matrix
of first partial derivatives, with the off-diagonal elements
multiplied by a factor of one half.

Dawid and Sebastiani (1999) and Pukelsheim (2006) give various
examples of convex functions $\phi$ whose domain is the cone of the
symmetric and positive definite elements of $\real^{q \times q}$, with
the matrix norm
\begin{equation}  \label{eq:matrix.norm} 
\phi(\Sigma) = \left( \, \frac{1}{q} \, {\rm tr} (\Sigma^s) \right)^{\! 1/s} 
\end{equation} 
for $s > 1$ being one such instance.  The matrix norm is nonnegative,
nondecreasing in the Loewner order, continuous, strictly convex,
standardized and homogeneous of order one.  With simple adaptations,
the construction extends to any real or extended real-valued exponent
$s$ and to general, not necessarily positive definite symmetric
matrices (Pukelsheim 2006, pp.~141 and 151).  In the limit as $s \to
0$ in (\ref{eq:matrix.norm}) the log determinant $\phi(\Sigma) = \log
\hsp \det (\Sigma)$ emerges.  When used in the Bregman representation
(\ref{eq:Bregman.matrix}), the log determinant function gives rise to
a well known homogeneous scoring function for point predictions of a
positive definite symmetrically matrix-valued quantity in $\real^{q
\times q}$, namely,
\begin{equation}  \label{eq:JamesStein}
\myS(\Sigma_x,\Sigma_y) = {\rm tr} \left( \Sigma_x^{-1} \Sigma_y \right) 
  - \log \det \left( \Sigma_x^{-1} \Sigma_y \right) - q, 
\end{equation} 
which was introduced by James and Stein (1961, Section 5).  When $q =
1$ the scoring function (\ref{eq:JamesStein}) reduces to the Patton
function (\ref{eq:Patton}) with $b = 0$, that is, the QLIKE function.

In the case of quantiles, the passage from the univariate functional
to multivariate analogues is much less straightforward.  Notions of
quantiles for multivariate distributions based on loss or scoring
functions have been studied by Abdous and Theodorescu (1992),
Chaudhuri (1996), Koltchinskii (1997), Serfling (2002) and Hallin,
Paindaveine and \u{S}iman (2010), among others.  In particular, it is
customary to define the {\em median}\/ of a probability distribution
$F$ on $\real^d$ as
\[
\textstyle
\hat{x} = \arg \min_x \hsp \myE_F \! \left( \| x - Y \| - \| Y \| \right) \! ,    
\]
where $\| \cdot \|$ denotes the Euclidean norm (Small 1990).  If $d =
1$, this yields the traditional median on the real line, with the $\|
Y \|$ term eliminating the need for moment conditions on the
predictive distribution (Kemperman 1987).  Of course, norms and
distances other than the Euclidean could be considered.  In this more
general type of situation, Koenker (2006) proposed that a functional
based on minimizing the square of a distance be called a {\em
Fr\'echet mean}, and a functional based on minimizing a distance a
{\em Fr\'echet median}, just as in the traditional case of the
Euclidean distance.

\section{Discussion} \label{sec:discussion} 

Ideally, forecasts ought to be probabilistic, taking the form of
predictive distributions over future quantities and events (Dawid
1984; Diebold et al.~1998; Granger and Pesaran 2000a, 2000b; Gneiting
2008a).  If point forecasts are to be issued and evaluated, it is
essential that either the scoring function be specified ex ante, or an
elicitable target functional be named, such as the mean or a quantile
of the predictive distribution, and scoring functions be used that are
consistent for the target functional.

Our plea for the use of consistent scoring functions supplements and
qualifies, but does not contradict, extant recommendations in the
forecasting literature, such as those of Armstrong (2001), Jolliffe
and Stephenson (2003) and Fildes and Goodwin (2007).  For example,
Fildes and Goodwin (2007) propose forecasting principles for
organizations, the eleventh of which suggests that ``multiple measures
of forecast accuracy'' be employed.  I agree, with the qualification
that the scoring functions to be used be consistent for the target
functional.

We have developed theory for the notions of consistency and
elicitability, and have characterized the classes of the loss or
scoring functions that result in expectations, ratios of expectations,
quantiles or expectiles as optimal point forecasts.  Some of these
results are classical, such as those for means and quantiles (Savage
1971; Thomson 1979), while others are original, including a
disconcerting negative result, in that scoring functions which are
consistent for the CVaR functional do not exist.

In the case of the mean functional, the consistent scoring functions
are the Bregman functions of the form (\ref{eq:Bregman}).  Among
these, a particularly attractive choice is the Patton family
(\ref{eq:Patton}) of homogeneous scoring functions, which nests the
squared error (SE) and QLIKE functions.  In evaluating volatility
forecasts, Patton and Sheppard (2009) recommend the use of the latter
because of its superior power in Diebold and Mariano (1995) and West
(1996) tests of predictive ability, which depend on differences
between mean scores of the form (\ref{eq:S}) as test statistics.
Further work in this direction is desirable, both empirically and
theoretically.  If quantile forecasts are to be assessed, the
consistent scoring functions are the GPL functions of the form
(\ref{eq:GPL}), with the homogeneous power functions in
(\ref{eq:GPL.b}) being appealing examples.  Interestingly, the scoring
functions that are consistent for expectiles combine key elements of
the Bregman and GPL families.

As regards the most commonly used scoring functions in academia,
businesses and organizations, the squared error scoring function is
consistent for the mean, and the absolute error scoring function for
the median.  The absolute percentage error scoring function, which is
commonly used by businesses and organizations, and occasionally in
academia, is consistent for a non-standard functional, namely, the
median of order $-1$, ${\rm med}^{(-1)}$, which tends to support
severe underforecasts, as compared to the mean or median.  It thus
seems prudent that businesses and organizations consider the intended
or unintended consequences and reassess its suitability as a scoring
function.

Pers et al.~(2009) propose a game of prediction for a fair comparison
between competing predictive models, which employs proper scoring
rules.  As Theorem \ref{th:proper} shows, consistent scoring functions
can be interpreted as proper scoring rules.  Hence, the protocol of
Pers et al.~(2009) applies directly to the evaluation of point
forecasting methods.  Their focus is on the comparison of custom-built
predictive models for a specific purpose, as opposed to the
M-competitions in the forecasting literature (Makridakis and Hibon
1979, 2000; Makridakis et al.~1982, 1993), which compare the
predictive performance of point forecasting methods across multiple,
unrelated time series.  In this latter context, additional
considerations arise, such as the comparability of scores across time
series with realizations of differing magnitude and volatility, and
commonly used evaluation methods remains controversial (Armstrong and
Collopy 1992; Fildes 1992; Ahlburg et al.~1992; Hyndman and Koehler
2006).

The notions of consistency and elicitability apply to point forecast
competitions, where participants ought to be advised ex ante about
the scoring function(s) to be employed, or, alternatively, target
functional(s) ought to be named.  If multiple target functionals are
named, participants can enter possibly distinct point forecasts for
distinct functionals.  Similarly, if multiple scoring functions are to
be used in the evaluation, and the scoring functions are consistent
for distinct functionals, participants ought to be allowed to submit
possibly distinct point forecasts.

While thus far we have addressed forecasting or prediction problems,
similar issues arise when the goal is estimation.  Technically, our
discussion relates to M-estimation (Huber 1964; Huber and Ronchetti
2009).  A century ago Keynes (1911, p.~325) derived the Bregman
representation (\ref{eq:Bregman}) in characterizing the probability
density functions for which the ``most probable value'' is the
arithmetic mean.  For a contemporary perspective in terms of maximum
likelihood and M-estimation, see Klein and Grottke (2008).  Komunjer
(2005) applied the GPL class (\ref{eq:GPL}) in conditional quantile
estimation, in generalization of the traditional approach to quantile
regression, which is based on the asymmetric piecewise linear scoring
function (Koenker and Bassett 1978).  Similarly, Bregman functions of
the original form (\ref{eq:Bregman}) and of the variant in
(\ref{eq:Bregman.tau}) could be employed in generalizing symmetric and
asymmetric least squares regression.

In applied settings, the distinction between prediction and estimation
is frequently blurred.  For example, Shipp and Cohen (2009) report on
U.S.~Census Bureau plans for evaluating population estimates against
the results of the 2010 Census.  Five measures of accuracy are to be
used to assess the Census Bureau estimates, including the root mean
squared error (SE) and the mean absolute percentage error (APE).  Our
results demonstrate that Census Bureau scientists face an impossible
task in designing procedures and point estimates aimed at minimizing
both measures simultaneously, because the SE and the APE are
consistent for distinct statistical functionals.  In this light, it
may be desirable for administrative or political leadership to provide
a directive or target functional to Census Bureau scientists, much in
the way that Murphy and Daan (1985) and Engelberg et al.~(2009)
requested guidance for point forecasters, in the quotes that open and
motivate this paper.

\section*{Appendix A: Proofs}

\smallskip
{\em Proof of Theorem \ref{th:convex}.}  Given $F \in \cF$, let $t \in
\myT(F)$ and $x \in \myD$.  Then
\begin{eqnarray*} 
\myE_F \, \myS(t,Y) 
& = & \myE_F \! \left[ \, \int \myS_{\hsp \omega}(t,y) \: \lambda(\dd \omega) \right] \\
& = & \int \left[ \myE_F \rule{0mm}{4mm} \, \myS_{\hsp \omega}(t,y) \right] \lambda(\dd \omega) \\
& \leq & \int \left[ \myE_F \rule{0mm}{4mm} \, \myS_{\hsp \omega}(x,y) \right] \lambda(\dd \omega)
\;\; = \;\; \myE_F \, \myS(x,Y),  
\end{eqnarray*} 
where the interchange of the expectation and the integration is
allowable, because each $\myS_\omega$ is a nonnegative scoring
function.  \done

\bigskip
{\em Proof of Theorem \ref{th:proper}.}  Given any two probability
measures $F, G \in \cF$, we have
\begin{eqnarray*}
\myE_F \, \bS(F,Y) = \myE_F \, \myS(\myT(F),Y) 
\leq \myE_F \, \myS(\myT(G),Y) = \myE_F \, \bS(G,Y),   
\end{eqnarray*}  
where the expectations are well-defined, because the scoring function
$\myS$ is nonnegative.  \done

\bigskip
{\em Proof of Theorem \ref{th:relevation}.}  We first show part (b).
Towards this end, let $t_g \in \myT_g(F)$ and $x_g \in \myD$.  Then
$t_g = g(t)$ for some $t \in \myT(F)$ and $x_g = g(x)$ for some $x \in
\myD$.  Therefore,
\[
\myE_F \, \myS_g(t_g,Y) 
= \myE_F \, \myS( \hsp t,Y)
\leq \myE_F \, \myS(x,Y)
= \myE_F \, \myS_g(x_g,Y).   
\]
As regards parts (c) and (a), it suffices to note that if $\myS$ is
strictly consistent, we have equality if and only if $x \in \myT(F)$
or, equivalently, $x_g \in \myT_g(F)$.  \done

\bigskip
{\em Proof of Theorem \ref{th:weighted}.}  We first prove part (b).
Let $F \in \cF^{(w)}$, $t \in \myT^{(w)}(F)$ and $x \in \myD$.  Then
\begin{eqnarray*} 
\myE_F \, \myS^{(w)}(t,Y) 
& = & \myE_F \! \hsp \left[ w(Y) \, \myS(t,Y) \right] \\  
& = & \int \myS(t,y) \hsp w(y) \hsp f(y) \: \mu(\dd y) \\
& = & \left[ \, \int \myS(t,y) \, \dd F^{(w)}(y) \right] \cdot
      \left[ \, \int w(y) \hsp f(y) \: \mu(\dd y) \right]^{-1} \\    
& \leq & \left[ \, \int \myS(x,y) \, \dd F^{(w)}(y) \right] \cdot 
         \left[ \, \int w(y) \hsp f(y) \: \mu(\dd y) \right]^{-1} \\
& = & \myE_F \! \hsp \left[ w(Y) \, \myS(x,Y) \right] \rule{0mm}{5.5mm} \\       
& = & \myE_F \! \hsp \left[ \myS^{(w)}(x,Y) \right] \! , \rule{0mm}{6.5mm}       
\end{eqnarray*} 
where $\mu$ is a dominating measure.  The critical inequality holds
because $F^{(w)} \in \cF^{(w)} \subseteq \cF$ and $t^{(w)} \in
\myT^{(w)}(F) = \myT(F^{(w)})$.  To prove parts (c) and (a), we note
that the inequality is strict if $\myS$ is strictly consistent for
$\myS$, unless $x \in \myT(F^{(w)}) = \myT^{(w)}(F)$.  \done

\bigskip
{\em Proof of Theorem \ref{th:necessary}.}  Suppose that the
functional $\myT$ is elicitable relative to the class $\cF$ on the
domain $\myD$.  Then there exists a scoring function $\myS$ which is
strictly consistent for it relative to $\cF$.  Suppose now that $F_0
\in \cF$, $F_1 \in \cF$ and $t \in \myD$ are such that $t \in
\myT(F_0)$ and $t \in \myT(F_1)$.  If $x \in \myD$ is arbitrary and $p
\in (0,1)$ is such that $F_p = (1-p)F_0 + p F_1 \in \cF$ then
\begin{eqnarray*} 
\myE_{F_p} \hsp \myS(t,Y) & = & (1-p) \, \myE_{F_0} \myS(t,Y) + p \, \myE_{F_1} \myS(t,Y) \\
& \leq & (1-p) \, \myE_{F_0} \myS(x,Y) + p \, \myE_{F_1} \myS(x,Y) 
\;\; = \;\; \myE_{F_p} \hsp \myS(x,Y).  
\end{eqnarray*} 
Hence, $t \in \myT(F_p)$. 
\done

\bigskip
{\em Sketch of the proof of Theorem \ref{th:mean}.}  The statements in
parts (b) and (c) are immediate from the arguments in Section 6.3 of
Savage (1971), and form special cases of the more general result in
Theorem \ref{th:ratio}.  To prove the necessity of the representation
(\ref{eq:Bregman}), Savage essentially applied Osband's principle with
the identification function $\myV(x,y) = x - y$. \done

\bigskip
{\em Proof of Theorem \ref{th:ratio}.}  We first prove part (b).  To
show the sufficiency of the representation (\ref{eq:Bregman.ratio}),
let $x \in \myI$ and let $F$ be a probability measure on $\myI$ for
which $\myE_F \hsp [ \hsp r(Y)]$, $\myE_F \hsp [ \hsp s(Y)]$, $\myE_F
\hsp [ \hsp r(Y) \hsp \phi'(Y)]$, $\myE_F \hsp [ \hsp s(Y) \hsp
\phi(Y)]$ and\/ $\myE_F \hsp [ \hsp\hsp Y s(Y) \hsp \phi'(Y)]$ exist
and are finite.  Then
\begin{eqnarray*} 
\lefteqn{\hspace{-2.5mm} \myE_F \, \myS(x,Y) - \myE_F \, \myS \! \left(
\Ers, Y \right)} \\ & = & \myE_F \hsp [s(Y)] \left\{ \phi \! \left(
\Ers \right) - \phi(x) - \phi'(x) \left( \Ers - x \right) \right\}
\end{eqnarray*} 
is nonnegative, and is strictly positive if $\phi$ is strictly convex
and $x \not= \myE_F \hsp [r(Y)] \, / \, \myE_F \hsp [s(Y)]$.

As regards part (c), it remains to show the necessity of the
representation (\ref{eq:Bregman.ratio}).  We apply Osband's principle
with the identification function $\myV(x,y) = x \hsp s(y) - r(y)$, as
proposed by Osband (1985, p.~14).  Arguing in the same way as in
Section \ref{sec:Osband}, we see that
\[
\myS_{(1)}(x,a) \hsp / \hsp (x \hsp s(a) - r(a)) = 
\myS_{(1)}(x,b) \hsp / \hsp (x \hsp s(b) - r(b))
\]
for all pairwise distinct $a, b$ and $x \in \myI$.  Hence,
\[
\myS_{(1)}(x,y) = h(x) \hsp (x \hsp s(y) - r(y)) 
\]
for $x, y \in \myI$ and some function $h : \myI \to \myI$.  Partial
integration yields the representation (\ref{eq:Bregman.ratio}), where
\begin{equation}  \label{eq:phi} 
\phi(x) = \int_{x_0}^x \int_{x_0}^s h(u) \: \dd u \: \dd s
\end{equation} 
for some $x_0 \in \myI$.  Finally, $\phi$ is convex, because the
scoring function $\myS$ is nonnegative, which implies the validity of
the subgradient inequality.

To prove part (a), we consider the scoring function
(\ref{eq:Bregman.ratio}) with $\phi(y) = y^2/(1+|y|)$, for which the
expectations in part (b) exist and are finite if, and only if, $\myE_F
\hsp [ \hsp r(Y)]$, $\myE_F \hsp [ \hsp s(Y)]$ and\/ $\myE_F \hsp [
\hsp\hsp Y s(Y) \hsp ]$ exist and are finite. \done

\bigskip
{\em Sketch of the proof of Theorem \ref{th:quantiles}.}  For concise
yet full-fledged proofs of parts (b) and (c), see Gneiting (2008b),
where Osband's principle is applied with the identification function
$\myV(x,y) = \one(x \geq y) - \alpha$.  To prove part (a), we may
apply part (c) with any strictly increasing, bounded function $g :
\myI \to \myI$, with $g(x) = \exp(-x)/(1+\exp(-x))$ being one such
example. \done

\bigskip
{\em Proof of Theorem \ref{th:expectiles}.}  To show the sufficiency
of the representation (\ref{eq:Bregman.tau}), let $x \in \myI$ where
$x < \mu_\tau$, and let $F$ be a probability measure with compact
support in $\myI$.  A tedious but straightforward calculation shows
that if $\myS$ is of the form (\ref{eq:Bregman.tau}) then
\begin{eqnarray*} 
\lefteqn{\hspace{-2.5mm} \myE_F \, \myS(x,Y) - \myE_F \, \myS \! \left( \mu_\tau, Y \right)} \\ & = & 
(1 - \tau) \int_{(-\infty, \, x)} \left( \phi(\mu_\tau) - \phi(x) - \phi'(x) (\mu_\tau - x) \right) \dd F(y) \\
&& + \; \tau \int_{[x, \, \mu_\tau)} \left( \phi(y) - \phi(x) - \phi'(x) (y - x) \right) \dd F(y) \\
&& + \; \tau \int_{[\mu_\tau, \, \infty)} \left( \phi(\mu_\tau) - \phi(x) - \phi'(x) (\mu_\tau - x) \right) \dd F(y) \\
&& + \; (1 - \tau) \int_{[x, \, \mu_\tau)}  
\underbrace{\left( \phi(\mu_\tau) - \phi(y) - \phi'(x) (\mu_\tau - y) \right)}
_{\geq \; \phi(\mu_\tau) \, - \, \phi(y) \, - \, \phi'(y) \hsp (\mu_\tau - y) \; \geq \; 0} \dd F(y) 
\end{eqnarray*} 
is nonnegative, and is strictly positive if $\phi$ is strictly convex.
An analogous argument applies when $x > \mu_\tau$.  This proves
sufficiency in part (b) as well as the claim in part (c).

To prove the necessity of the representation (\ref{eq:Bregman.tau}) in
part (b), we apply Osband's principle with the identification function
$\myV(x,y) = |\one(x \geq y) - \tau \hsp | \, (x - y)$.  Arguing in
the usual way, we see that
\[
\myS_{(1)}(x,y) = h(x) \, V(x,y) 
\]
for $x, y \in \myI$ and some function $h : \myI \to \myI$.  Partial
integration yields the representation (\ref{eq:Bregman.tau}), where
$\phi$ is defined as in (\ref{eq:phi}) and is convex, because $\myS$
is nonnegative.

To prove part (a), we apply part (c) with the convex function $\phi(y)
= y^2/(1+|y|)$, for which $\myE_F \, \phi(Y)$ exists and is finite if,
and only if, $\myE_F \hsp Y$ exists and is finite.  \done

\bigskip
{\em Proof of Theorem \ref{th:CVaR}.}  Suppose first that $\cF$
contains the measures with finite support.  Let $a, b, c, d \in \myI$
be such that $a < b < c < \frac{1}{2}(b+d)$, which implies $b < d$,
and consider the probability measures
\[
F_1 = \alpha \hsp\hsp \delta_a + \frac{1}{2} \hsp (1 - \alpha) \left( \delta_{b} + \delta_{d} \right) \! ,  
\qquad
F_2 = \alpha \hsp\hsp \delta_c  + (1 - \alpha) \hsp\hsp \delta_{(b+d)/2},  
\]
where $\delta_x$ denotes the point measure in $x \in \real$.  Then
$\mbox{CVaR}_\alpha(F_1) = \mbox{CVaR}_\alpha(F_2) =
\frac{1}{2}(b+d)$, while $\mbox{CVaR}_\alpha(\frac{1}{2}(F_1+F_2)) =
\frac{1}{4}(b+c+2d) > \frac{1}{2}(b+d)$.  Thus, the level sets of the
functional are not convex.  By Theorem \ref{th:necessary}, the CVaR
functional is not elicitable relative to the class $\cF$.  An
analogous example emerges when the point measures are replaced by
appropriately focused and centered absolutely continuous distributions
with compact support.  \done

\section*{Appendix B: Optimal point forecasts under the relative error scoring function (Table \ref{tab:Patton})}

Here we address a problem posited by Patton (2010), in that we find
the optimal point forecast or Bayes rule
\begin{equation}  \label{eq:Bayes.RE} 
\textstyle \hat{x} = \arg \min_{\hsp x} \hsp 
\displaystyle \myE_F \, \myS(x,y)
\quad \mbox{under} \quad \myS(x,y) = |(x-y)/x|, 
\end{equation} 
where $Y = Z^2$ and $Z$ has a $t$-distribution with mean 0, variance 1
and $\nu > 2$ degrees of freedom.  In the limiting case as $\nu \to
\infty$, we take $Z$ to be standard normal.

To find the optimal point forecast, we apply Theorem \ref{th:dual} and
part (b) of Theorem \ref{th:weighted} with the original scoring
function $\myS(x,y) = |x^{-1} - y^{-1}|$, the weight function $w(y) =
y$ and the domain $\myD = (0,\infty)$, so that $\myS^{(w)}(x,y) =
|(x-y)/x|$.  By Theorem \ref{th:quantiles}, the scoring function
$\myS$ is consistent for the median functional.  Therefore, by Theorem
\ref{th:weighted} the optimal point forecast under the weighted
scoring function $\myS^{(w)}$ is the median of the probability
distribution whose density is proportional to $y \hsp f(y)$, where $f$
is the density of $Y$, or equivalently, proportional to $y^{1/2} \hsp
g(y^{1/2})$, where $g$ is the density of $Z$.

Hence, if $Z$ has a $t$-distribution with mean 0, variance 1 and $\nu
> 2$ degrees of freedom, the optimal point forecast under the relative
error scoring function is the median of the probability distribution
whose density is proportional to
\[
y^{1/2} \left( 1 + \frac{y}{\nu-2} \right)^{- \, (\nu+1)/ \, 2}     
\] 
on the positive halfaxis.  Using any computer algebra system, this
median can readily be computed symbolically or numerically, to any
desired degree of accuracy.  For example, if $\nu = 4$ the optimal
point forecast (\ref{eq:Bayes.RE}) is
\[
\hat{x} = \frac{2}{2^{2/3}-1} = 3.4048 \ldots 
\]
Table \ref{tab:Patton} provides numerical values along with the
approximations in Table 1 of Patton (2010), which were obtained by
Monte Carlo methods, and thus are less accurate.  If $Z$ has variance
$\sigma^2$, the entries in the table continue to apply, if they are
multiplied by this constant.

\section*{Acknowledgements}

The author thanks Werner Ehm, Marc G.~Genton, Peter Guttorp, Jorgen
Hilden, Peter J.~Huber, Ian T.~Jolliffe, Charles F.~Manski, Caren
Marzban, Kent H.~Osband, Pierre Pinson, Adrian E.~Raftery, Ken Rice,
R.~Tyrrell Rockafellar, Paul D.~Sampson, J.~McLean Sloughter, Stephen
Stigler, Adam Szpiro, Jon A.~Wellner and Robert L.~Winkler for
discussions, references and preprints.  Financial support was provided
by the Alfried Krupp von Bohlen und Halbach Foundation, and by the
National Science Foundation under Awards ATM-0724721 and DMS-0706745
to the University of Washington.  Special thanks go to University of
Washington librarians Martha Tucker and Saundra Martin for their
unfailing support of the literature survey in Table \ref{tab:lit}.  Of
course, the opinions expressed in this paper as well as any errors are
solely the responsibility of the author.

\section*{References}

\newenvironment{reflist}{\begin{list}{}{\itemsep 0mm \parsep 1mm
\listparindent -7mm \leftmargin 7mm} \item \ }{\end{list}}

\vspace{-6.5mm}
\begin{reflist}

{\rm Abdous, B., and Theodorescu, R.}~(1992), ``Note on the Spatial
Quantile of a Random Vector,'' {\em Statistics \& Probability
Letters}, 13, 333--336.

{\rm Acerbi, C.}~(2002), ``Spectral Measures of Risk: A Coherent
Representation of Subjective Risk Aversion,'' {\em Journal of Banking
and Finance}, 26, 1505--1518.

{\rm Ahlburg, D.~A., Chatfield, C., Taylor, S.~J., Thompson, P.~H.,
Murphy, A.~H., Winkler, R.~L., Collopy, F., Armstrong, J.~S.~and
Fildes, R.}~(1992), ``A Commentary on Error Measures,'' {\em
International Journal of Forecasting}, 8, 99--111.

{\rm Armstrong, J.~S.}~(2001), ``Evaluating Forecasting Methods,'' in
{\em Principles of Forecasting}, Armstrong, J.~S., ed., Kluwer, 
Norwell, Massachusetts, pp.~443--471.

{\rm Armstrong, J.~S., and Collopy, F.}~(1992), ``Error Measures for
Generalizing About Forecasting Methods: Empirical Comparisons,'' {\em
International Journal of Forecasting}, 8, 69--80.

{\rm Artzner, P., Delbaen, F., Eber, J.-M.~and Heath, D.}~(1999),
``Coherent Measures of Risk,'' {\em Mathematical Finance}, 9,
203--228.

{\rm Banerjee, A., Guo, X.~and Wang, H.}~(2005), ``On the Optimality of
Conditional Expectation as a Bregman Predictor,'' {\em IEEE
Transactions on Information Theory}, 51, 2664--2669.

{\rm Bao, Y., Lee, T.-H., and Salto\u{g}lu, B.}~(2006), ``Evaluating
Predictive Performance of Value-at-Risk Models in Emerging Markets: A
Reality Check,'' {\em Journal of Forecasting}, 25, 101--128.

{\rm Berkowitz, J., and O'Brien, J.}~(2002), ``How Accurate are
Value-at-Risk Models at Commercial Banks?,'' {\em Journal of Finance},
57, 1093--1111.

{\rm Bollerslev, T.}~(1986), ``Generalized Autoregressive Conditional
Heteroscedasticity,'' {\em Journal of Econometrics}, 31, 307--327.

{\rm Buja, A., Stuetzle, W.~and Shen, Y.}~(2005), ``Loss Functions for
Binary Class Probability Estimation and Classification: Structure and
Applications,'' Working paper, \newline
\url{http://www-stat.wharton.upenn.edu/~buja/PAPERS/paper-proper-scoring.pdf}.

{\rm Carbone, R., and Armstrong, J.~S.}~(1982), ``Evaluation of
Extrapolative Forecasting Methods: Results of a Survey of Academicians
and Practicioners,'' {\em Journal of Forecasting}, 1, 215--217.

{\rm Cervera, J.~L., and Mu\~{n}oz, J.}~(1996), ``Proper Scoring Rules
for Fractiles,'' in {\em Bayesian Statistics 5}, Bernardo, J.~M.,
Berger, J.~O., Dawid, A.~P., and Smith, A.~F.~M., eds., Oxford
University Press, pp.~513--519.

{\rm Chaudhuri, P.}~(1996), ``On a Geometric Notion of Quantiles for
Multivariate Data,'' {\em Journal of the American Statistical
Association}, 91, 862--872.

{\rm Christoffersen, P.~F., and Diebold, F.~X.}~(1996), ``Further
Results on Forecasting and Model Selection Under Asymmetric Loss,''
{\em Journal of Applied Econometrics}, 11, 561--571.

{\rm Dawid, A.~P.}~(1984), ``Statistical Theory: The Prequential
Approach,'' {\em Journal of the Royal Statistical Society, Ser.~A},
147, 278--292.

{\rm \bysame}~(2007), ``The Geometry of Proper Scoring Rules,''
{\em Annals of the Institute of Statistical Mathematics}, 59, 77--93.

{\rm Dawid, A.~P.~and Sebastiani, P.}~(1999), ``Coherent Dispersion
Criteria for Optimal Experimental Design,'' {\em Annals of
Statistics}, 27, 65--81.

{\rm DeGroot, M.~H., and Fienberg, S.~E.}~(1983), ``The Comparison and
Evaluation of Probability Forecasters,'' {\em Statistician}, 12,
12--22.

{\rm Diebold, F.~X., and Mariano, R.~S.}~(1995), ``Comparing Predictive
Accuracy,'' {\em Journal of Business and Economic Statistics}, 13,
253--263.

{\rm Diebold, F.~X., Gunther, T.~A., and Tay, A.~S.}~(1998), ``Evaluating
Density Forecasts With Applications to Financial Risk Management,''
{\em International Economic Review}, 39, 863--883.

{\rm Duffie, D., and Pan, J.}~(1997), ``An Overview of Value at Risk,''
{\em Journal of Derivatives}, 4, 7--49.

{\rm Efron, B.}~(1991), ``Regression Percentiles Using Asymmetric
Squared Error Loss,'' {\em Statistica Sinica}, 1, 93--125.

{\rm Engelberg, J., Manski, C.~F., and Williams, J.}~(2009),
``Comparing the Point Predictions and Subjective Probability
Distributions of Professional Forecasters,'' {\em Journal of Business
and Economic Statistics}, 27, 30--41.

{\rm Engle, R.~F.}~(1982), ``Autoregressive Conditional
Heteroscedasticity With Estimates of the Variance of United Kingdom
Inflation,'' {\em Econometrica}, 45, 987--1007.

{\rm Ferguson, T.~S.}~(1967), {\em Mathematical Statistics: A
Decision-Theoretic Approach}, Academic, New York.

{\rm Fildes, R.}~(1992), ``The Evaluation of Extrapolative Forecasting
Methods,'' {\em International Journal of Forecasting}, 8, 81--98.

{\rm Fildes, R., and Goodwin, P.}~(2007), ``Against Your Better
Judgement? How Organizations Can Improve Their Use of Management
Judgement in Forecasting,'' {\em Interfaces}, 37, 570--576.

{\rm Fildes, R., Nikolopoulos, K., Crone, S.~F., and Syntetos,
A.~A.}~(2008), ``Forecasting and Operational Research: A Review,''
{\em Journal of the Operational Research Society}, 59, 1150--1172.

{\rm Giacomini, R., and Komunjer, I.}~(2005), ``Evaluation and
Combination of Conditional Quantile Forecasts,'' {\em Journal of
Business and Economic Statistics}, 23, 416--431.

{\rm Gneiting, T.}~(2008a), ``Editorial: Probabilistic Forecasting,''
{\em Journal of the Royal Statistical Society, Ser.~A}, 171, 319--321.

{\rm \bysame}~(2008b), ``Quantiles as Optimal Point Forecasts,''
Technical Report no.~538, University of Washington, Department of
Statistics, \newline
\url{http://www.stat.washington.edu/research/reports/2008/tr538.pdf}.

{\rm \bysame}~(2010), ``Quantiles as Optimal Point Forecasts,''
{\em International Journal of Forecasting}, in press.

{\rm Gneiting, T., and Raftery, A.~E.}~(2007), ``Strictly Proper
Scoring Rules, Prediction, and Estimation,'' {\em Journal of the
American Statistical Association}, 102, 359--378.

{\rm Gneiting, T., Stanberry, L.~I., Grimit, E.~P., Held, L., and
Johnson, N.~A.}~(2008), ``Assessing Probabilistic Forecasts of
Multivariate Quantities, With Applications to Ensemble Predictions of
Surface Winds,'' {\em Test}, 17, 211--264.

{\rm Granger, C.~W.~J.}~(1969), ``Prediction With a Generalized Cost
of Error Function'', {\em Operational Research Quarterly}, 20,
199--207.

{\rm Granger, C.~W.~J., and Pesaran, M.~H.}~(2000a), ``Economic and
Statistical Measures of Forecast Accuracy,'' {\em Journal of
Forecasting}, 19, 537--560.

{\rm \bysame}~(2000b), ``A Decision Theoretic Approach to Forecast
Evaluation,'' in {\em Statistics and Finance: An Interface}, Chan,
W.-S., Li, W.~K., and Tong, H., eds., Imperial College Press, London,
pp.~261--278.

{\rm Hallin, M., Paindaveine, D., and \u{S}iman, M.}~(2010),
``Regression Quantiles: From $L_1$ Optimization to Halfspace Depth,''
{\em Annals of Statistics}, 38, 635--703.

{\rm Henderson, R., Jones, M., and Stare, J.}~(2001), ``Accuracy of
Point Predictions in Survival Analysis,'' {\em Statistics in
Medicine}, 20, 3083--3096.

{\rm Hering, A.~S., and Genton, M.~G.}~(2010), ``Powering up with
Space-Time Wind Forecasting,'' {\em Journal of the American
Statistical Association}, in press.

{\rm Hilden, J.}~(2008), ``Scoring Rules for Evaluation of
Prognosticians and Prognostic Rules,'' Workshop notes, 
\url{http://biostat.ku.dk/~jh/}.

{\rm Horowitz, J.~L., and Manski, C.~F.}~(2006): ``Identification and
Estimation of Statistical Functionals Using Incomplete Data,'' {\em
Journal of Econometrics}, 132, 445--459.

{\rm Huber, P.~J.}~(1964), ``Robust Estimation of a Location Parameter,''
{\em Annals of Mathematical Statistics}, 35, 73--101.

{\rm Huber, P.~J., and Ronchetti, P.~M.}~(2009), {\em Robust
Statistics}, 2nd edition, Wiley, Hoboken, New Jersey.

{\rm Hyndman, R.~J., and Koehler, A.~B.}~(2006), ``Another Look at
Measures of Forecast Accuracy,'' {\em International Journal of
Forecasting}, 22, 679--688.

{\rm James, W., and Stein, C.}~(1961), ``Estimation With Quadratic
Loss,'' in {\em Proceedings of the Fourth Berkeley Symposium on
Mathematical Statistics and Probability, Vol.~1}, Neyman, J., ed.,
University of California Press, pp.~361--379.

{\rm Jolliffe, I.~T.}~(2008), ``The Impenetrable Hedge: A Note on
Propriety, Equitability and Consistency,'' {\em Meteorological
Applications}, 15, 25--29.

{\rm Jolliffe, I.~T., and Stephenson, D.~B., eds.}~(2003), {\em
Forecast Verification: A Practicioner's Guide in Atmospheric Science},
Wiley, Chichester.

{\rm Jose, V.~R.~R., and Winkler, R.~L.}~(2009), ``Evaluating Quantile
Assessments,'' {\em Operations Research}, 57, 1287--1297. 

{\rm Kemperman, J.~H.~B.}~(1987), ``The Median of a Finite Measure on
a Banach Space,'' in {\em Statistical Data Analysis Based on the $L_1$
Norm and Related Methods}, Dodge, Y., ed., North Holland,
pp.~217--230.

{\rm Keynes, J.~M.}~(1911), ``The Principal Averages and the Laws of
Error which Lead to Them,'' {\em Journal of the Royal Statistical
Society}, 74, 322--331.

{\rm Klein, I.~and Grottke, M.}~(2008), ``On J.~M.~Keynes' ``The
Principal Averages and the Laws of Error which Lead to Them'' --
Refinement and Generalisation.''  Discussion Paper, \newline 
\url{http://www.iwqw.wiso.uni-erlangen.de/forschung/07-2008.pdf}.

{\rm Koenker, R.}~(2005), {\em Quantile Regression}, Cambridge
University Press.

{\rm \bysame (2006)}, ``The Median is the Message: Toward the
Fr\'echet Mean,'' {\em Journal de la Soci\'et\'e Fran\c{c}aise de
Statistique}, 147, 61--64.

{\rm Koenker, R., and Bassett, G.}~(1978), ``Regression Quantiles,'' {\em
Econometrica}, 46, 33--50.

{\rm Koltchinskii, V.~I.}~(1997), ``M-Estimation, Convexity and
Quantiles,'' {\em Annals of Statistics}, 25, 435--477.

{\rm Komunjer, I.}~(2005), ``Quasi Maximum-Likelihood Estimation for 
Conditional Quantiles,'' {\em Journal of Econometrics}, 128, 137--164.

{\rm Lambert, N.~S., Pennock, D.~M., and Shoham, Y.}~(2008), ``Eliciting
Properties of Probability Distributions,'' Extended abstract,
Proceedings of the 9th ACM Conference on Electronic Commerce, July
8--12, 2008, Chicago, Illinois.

{\rm Laurent, S., Rombouts, J.~V.~K., and Violante, F.}~(2009), ``On
Loss Functions and Ranking Forecasting Performances of Multivariate
Volatility Models'', Discussion Paper, \newline
\url{http://www.cirpee.org/fileadmin/documents/Cahiers_2009/CIRPEE09-48.pdf}.

{\rm Lehmann, E.~L.}~(1951), ``A General Concept of Unbiasedness,''
{\em Annals of Mathematical Statistics}, 22, 587--592.

{\rm Lehmann, E., and Casella, G.}~(1998), {\em Theory of Point
Estimation}, 2nd edition, Springer, New York.

{\rm Makridakis, S., and Hibon, M.}~(1979), ``Accuracy of Forecasting:
An Empirical Investigation'' (with discussion), {\em Journal of the
Royal Statistical Society, Ser.~A}, 142, 97--145.

{\rm \bysame}~(2000), ``The M3-Competition: Results, Conclusions and
Implications,'' {\em International Journal of Forecasting}, 16,
451--476.

{\rm Makridakis, S., Chatfield, C., Hibon, M., Lawrance, M., Mills,
T., Ord, K., and Simmons, L.~F.}~(1993), ``The M2-Competition: A
Real-Time Judgementally Based Forecasting Study,'' {\em International
Journal of Forecasting}, 9, 5--22.

{\rm Makridakis, S., Andersen, A., Carbone, R., Fildes, R., Hibon, M.,
Lewan\-dowski, R., Newton, J., Parzen, E., and Winkler, R.}~(1982),
``The Accuracy of Extrapolation (Time Series) Methods: Results of a
Forecasting Competition,'' {\em Journal of Forecasting}, 1, 111--153.

{\rm McCarthy, J.}~(1956), ``Measures of the Value of Information,'' {\em
Proceedings of the National Academy of Sciences}, 42, 654--655.

{\rm McCarthy, T.~M., Davis, D.~F., Golicic, S.~L., and Mentzner,
  J.~T.}  (2006), ``The Evolution of Sales Forecasting Management: A
20-Year Longitudinal Study of Forecasting Practice,'' {\em Journal of
  Forecasting}, 25, 303--324.

{\rm Mentzner, J.~T., and Kahn, K.~B.}~(1995), ``Forecasting Technique
Familiarity, Satisfaction, Usage, and Application,'' {\em Journal of
Forecasting}, 14, 465--476.

{\rm Moskaitis, J.~R., and Hansen, J.~A.}~(2006), ``Deterministic
Forecasting and Verification: A Busted System?,'' Working paper,
Massachusetts Institute of Technology, \newline 
\url{http://wind.mit.edu/~hansen/papers/MoskaitisHansenWAF2006.pdf}.

{\rm Murphy, A.~H., and Daan, H.}~(1985), ``Forecast Evaluation,'' in
{\em Probability, Statistics and Decision Making in the Atmospheric
Sciences}, Murphy, A.~H., and Katz, R.~W., eds., Westview Press,
Boulder, Colorado, pp.~379--437.

{\rm Murphy, A.~H., and Winkler, R.~L.}~(1987), ``A General Framework
for Forecast Verification'', {\em Monthly Weather Review}, 115,
1330--1338.

{\rm Noorbaloochi, S., and Meeden, G.}~(1983), ``Unbiasedness as the
Dual of Being Bayes,'' {\em Journal of the American Statistical
Association}, 78, 619--623.

{\rm Newey, W.~K., and Powell, J.~L.}~(1987), ``Asymmetric Least
Squares Estimation and Testing,'' {\em Econometrica}, 55, 819--847.

{\rm Offerman, T., Sonnemans, J., van de Kuilen, G., and Wakker,
P.~P.}~(2009), ``A Truth-Serum for non-Bayesians: Correcting Proper
Scoring Rules for Risk Attitudes.  {\em Review of Economic Studies},
76, 1461--1489.

{\rm Osband, K.~H.}~(1985), ``Providing Incentives for Better Cost
Forecasting,'' Ph.D.~Thesis, University of California, Berkeley.

{\rm Osband, K., and Reichelstein, S.}~(1985), ``Information-Eliciting
Compensation Schemes,'' {\em Journal of Public Economics}, 27,
107--115.

{\rm Park, H., and Stefanski, L.~A.}~(1998), ``Relative-Error
Prediction,'' {\em Statistics \& Probability Letters}, 40, 227--236.

{\rm Patton, A.~J.}~(2010), ``Volatility Forecast Comparison Using
Imperfect Volatility Proxies,'' {\em Journal of Econometrics}, in
press, \url{http://econ.duke.edu/~ap172/}.

{\rm Patton, A.~J., and Sheppard, K.}~(2009), ``Evaluating Volatility
and Correlation Forecasts,'' in {\em Handbook of Financial Time
Series}, Anderson, T.~G., Davis, R.~A., Kreiss, J.-P., and Mikosch,
T., eds., Springer, pp.~801--838.

{\rm Pers, T.~H., Albrechtsen, A., Holst, C., S{\o}rensen, T.~I.~A., and
Gerds, T.~A.}~(2009), ``The Validation and Assessment of Machine
Learning: A Game of Prediction from High-Dimensional Data,'' {\em PLoS
ONE}, 4, e6287, {\tt doi:10.1371/journal.pone.0006287}.

{\rm Phelps, R.~R.}~(1966), {\em Lectures on Choquet's Theorem},
D.~Van Nostrand, Princeton.

{\rm Pukelsheim, F.}~(2006), {\em Optimal Design of Experiments}, SIAM
Classics edition, SIAM, Philadelphia.

{\rm Raiffa, H., and Schlaifer, R.}~(1961), {\em Applied Statistical
Decision Theory}, Colonial Press, Clinton.

{\rm Reichelstein, S., and Osband, K.}~(1984), ``Incentives in
Government Contracts,'' {\em Journal of Public Economics}, 24,
257--270.

{\rm Rockafellar, R.~T.}~(1970), {\em Convex Analysis}, Princeton
University Press.

{\rm Rockafellar, R.~T., and Uryasev, S.}~(2000), ``Optimization of
Conditional Value-at-Risk,'' {\em Journal of Risk}, 2, 21--42.

{\rm \bysame}~(2002), ``Conditional Value-at-Risk for General Loss
Distributions,'' {\em Journal of Banking and Finance}, 26, 1443--1471.

{\rm Saerens, M.}~(2000), ``Building Cost Functions Minimizing to Some
Summary Statistics,'' {\em IEEE Transactions on Neural Networks},
11, 1263--1271.

{\rm Savage, L.~J.}~(1971), ``Elicitation of Personal Probabilities
and Expectations,'' {\em Journal of the American Statistical
Association}, 66, 783--810.

{\rm Schervish, M.~J.}~(1989), ``A General Method for Comparing
Probability Assessors,'' {\em Annals of Statistics}, 17, 1856--1879.

{\rm Serfling, R.}~(2002), ``Quantile Functions for Multivariate
Analysis: Approaches and Applications,'' {\em Statistica Neerlandica},
56, 214--232.

{\rm Shipp, S., and Cohen, S.}~(2009), ``COPAFS Focuses on Statistical
Activities,'' {\em Amstat News}, August 2009, 15--18.

{\rm Small, C.~G.}~(1990), ``A Survey of Multidimensional Medians,''
{\em International Statistical Review}, 58, 263--277.

{\rm Thomson, W.}~(1979), ``Eliciting Production Possibilities From a
Well-Informed Manager,'' {\em Journal of Economic Theory}, 20,
360--380.

{\rm Wellner, J.~A.}~(2009), ``Statistical Functionals and the Delta
Method,'' Lecture no\-tes, \newline
\url{http://www.stat.washington.edu/people/jaw/COURSES/580s/581/LECTNOTES/ch7.pdf}.

{\rm West, K.~D.}~(1996), ``Asymptotic Inference About Predictive Ability,'' {\em
Econometrica}, 64, 1067--1084.

{\rm Winkler, R.~L.}~(1996), ``Scoring Rules and the Evaluation of
Probabilities'' (with discussion), {\em Test}, 5, 1--60.
\end{reflist}

\end{document}